\documentclass[a4paper,11pt]{article}
\usepackage{latexsym, amssymb, enumerate, amsmath,epsfig}
\usepackage[margin=1in]{geometry}
\usepackage{setspace,color}
\usepackage[titletoc]{appendix}
\usepackage{epstopdf}
\usepackage{graphicx}
\usepackage{pgfplots}
\usepackage{float}
\usepackage{caption}
\usepackage{subcaption}
\usepackage{grffile}
\usepackage{amsthm}
\usepackage[metapost]{mfpic}
\usepackage{booktabs}
\usepackage{verbatim}
\usepackage{lineno}
\newcommand{\x}{\mathbf{x}}

\newcommand{\n}{\mathbf{n}}

\theoremstyle{definition}

\theoremstyle{remark}
\newtheorem{rem}{Remark}

\begin{document}

\title{Solving Partial Differential Equations on Evolving Surfaces via the Constrained Least-Squares and Grid-Based Particle Method}

\author{
Ningchen Ying\thanks{Department of Mathematics, the Hong Kong University of Science and Technology, Clear Water Bay, Hong Kong. Email: {\bf nying@connect.ust.hk}}
\and
Shingyu Leung\thanks{Department of Mathematics, the Hong Kong University of Science and Technology, Clear Water Bay, Hong Kong. Email: {\bf masyleung@ust.hk}}
}
\date{}
\maketitle

\begin{abstract}
We present a framework for solving partial different equations on evolving surfaces. Based on the grid-based particle method (GBPM) \cite{leuzha0801}, the method can naturally resample the surface even under large deformation from the motion law. We introduce a new component in the local reconstruction step of the algorithm and demonstrate numerically that the modification can improve computational accuracy when a large curvature region is developed during evolution. The method also incorporates a recently developed constrained least-squares ghost sample points (CLS-GSP) formulation \cite{yin17,yinchuleu24}, which can lead to a better-conditioned discretized matrix for computing some surface differential operators. The proposed framework can incorporate many methods and link various approaches to the same problem. Several numerical experiments are carried out to show the accuracy and effectiveness of the proposed method.
\end{abstract}

\section{Introduction}

This paper proposes a new numerical method to solve a class of parabolic partial differential equations (PDEs) on evolving surfaces. In particular, we consider the following equation
\begin{eqnarray}
\frac{Du}{Dt} +u\nabla_{\Gamma} \cdot \mathbf{v} + \nabla_{\Gamma} \cdot \mathbf{q} = f  \text{\;  on\; } \Gamma \, ,
\label{ParabolicEqn}
\end{eqnarray}
where $\Gamma$ is a surface evolving under the velocity flow $\mathbf{v}$, $\frac{D}{Dt}$ is the material derivative given by $\frac{D}{Dt} = \frac{\partial}{\partial t} + \mathbf{v} \cdot \nabla$, the vector $\mathbf{q}$ represents a flux and $f$ is a source term. This equation can also be rewritten as

$$
u_t + V\frac{\partial u}{\partial \n} -\kappa V u + \nabla_{\Gamma}\cdot(u\mathbf{T}) + \nabla_{\Gamma} \cdot \mathbf{q} = f  \text{\;  on\; } \Gamma \, ,
$$
where the quantity $\kappa$ is the mean curvature of the surface, $V$ and $\mathbf{T}$ are the normal velocity of the interface motion and the velocity tangential to the interface, respectively, satisfying the decomposition $\mathbf{v} = V\n + \mathbf{T}$. For motions in the normal direction, the equation can be further reduced to

$$
u_t + V\frac{\partial u}{\partial \n} -\kappa V u + \nabla_{\Gamma} \cdot \mathbf{q} = f  \text{\;  on\; } \Gamma \, .
$$
We can use various choices of the flux function $\mathbf{q}$. Here we only consider two cases with the diffusive flux given by $\mathbf{q} = -\mathcal{D}\nabla_{\Gamma}w$ where $w$ is a function depending on the quantity $u$ defined on the interface and also a diffusion parameter $\mathcal{D}>0$. In particular, if we choose $w= u$, equation (\ref{ParabolicEqn}) leads to the advection-diffusion equation with a source term $f$ given by
\begin{eqnarray}
\frac{Du}{Dt} +u\nabla_{\Gamma} \cdot \mathbf{v} - \mathcal{D}\Delta_{\Gamma} u = f \, .
\label{Eqn.AdvectionDiffusion}
\end{eqnarray}
Another widely used choice of the flux is $w = -\mu\Delta_{\Gamma}u + \psi'(u)$ where $\psi(u)$ is the so-called double well potential function and $\mu>0$ represents physically a chemical potential in modeling multiphase fluid flows. This leads to the Cahn-Hilliard equation on moving surfaces

$$
\frac{Du}{Dt} +u\nabla_{\Gamma} \cdot \mathbf{v} - \mathcal{D}\Delta_{\Gamma} (-\mu\Delta_{\Gamma} u +\psi'(u)) = f \, .
$$

These equations have many applications and are especially important in biological, physical, and engineering sciences. For example, the growth of a solid tumor has been modeled by a coupled system of a reaction-diffusion system with a curvature-dependent evolution law 
\cite{cragafmai99,chagangra01,ellsty12}. The Cahn-Hilliard has been used for simulating pattern formation and selection \cite{madmai07,barellmad11}. The advection-diffusion equation can model the multi-phase interfacial flows with a surfactant as discussed in \cite{xuzha03,xllz06}. 

There are various difficulties in developing a numerical strategy for the problems. The representation of the interface itself already plays an essential role since the topology of the surface might change with evolution. A typical Lagrangian method based on a triangulation of the interface might require remeshing when the dynamic involves significant deformation. This might lead to some computational challenges in the implementation. Mesh surgery might also be needed to carefully reconnect the triangulation of the surface when there is any topological change in the evolution. Implicit methods based on the Eulerian formulation have recently become more attractive. One of the most widely used representations is based on the level set method \cite{oshset88,oshfed03}. The evolution of the interface is modeled by the numerical solution of a Hamilton-Jacobi equation, which is a nonlinear first-order PDE. More recently, we have developed a grid-based particle method (GBPM) and a cell-based particle method (CBPM) in \cite{leuzha0801,honleuzha14} to represent and model the interface dynamics. These grid-based or cell-based representations sample the surface by meshless Lagrangian particles according to an underlying uniform Eulerian mesh. This automatically provides a set of quasi-uniform sampling points of the interface. As the interface evolves according to the motion law, these unconnected non-parametrized sampling points will be continuously updated and can naturally capture changes in the interface topology.

The second phase in the development is to have a numerical approach for solving PDEs on a static surface given by a specific interface representation. For example, based on the implicit Eulerian representation, some projection operators have been introduced to extend a PDE off the interface to a small neighborhood \cite{bcos01,gre05}. Without the projection operator, the extension idea has also been used in \cite{wonleu15} for solving the surface eikonal equation using the fast sweeping method \cite{tcoz03,zha05}. Based on the level set representation, a so-called closest point method (CPM) has been developed in  \cite{ruumer08,macruu08,macruu09,macbraruu11} by replacing the surface PDE by an evolution equation approach. Some methods have been proposed for general surfaces represented by point clouds. For example, there are methods based on the radial basis functions formulation \cite{pir12,swkf15}, a local mesh method \cite{lailiazha13}, a virtual grid difference method \cite{wanleuzha17}, etc. 

A more challenging task is an approach to couple the above two tasks together based on techniques like operator-splitting so that one can develop a numerical method for solving PDEs on evolving surfaces. Some fully explicit Lagrangian methods have been developed in \cite{dziell07,barellmad11,dziell13} based on a finite element formulation on the triangulated surface. A fully implicit Eulerian method based on the level set method has been proposed in \cite{xuzha03}. A numerical approach based on the GBPM has been discussed in \cite{leuzha10,leulowzha11}. Coupling the GBPM, the CPM can also solve PDEs on moving surfaces \cite{petruu16}.

In this paper, we develop a general framework for solving PDEs on evolving surfaces represented by the GBPM representation. We closely follow the original GBPM in \cite{leuzha0801,leuzha10,leulowzha11} for both the interface representation and the evolution. Yet, we discuss an improvement in the implementation of the local reconstruction. We must collect a small set of neighbors in the original algorithm for the local interface reconstruction. These neighboring sample points must be chosen carefully for an accurate approximation. One constraint is that we must sample the same interface segment. When the curvature of the interface is too large compared to the underlying mesh size, one might collect footpoints from different segments of the interface, and it will result in a wrong local reconstruction of the surface. In the original work, we proposed checking if the associated normal vectors from these footpoints were pointing in a similar direction. When the interface developed a high curvature region compared to the size of the underlying mesh, however, this check does not provide enough penalty on far away sampling points. This paper introduces a simple weighting strategy that can naturally reduce the contributions from those footpoints from regions close to a different segment. 

The second main contribution of the paper is to incorporate a CLS-GSP approach for approximating local derivatives \cite{yin17,yinchuleu24}. In the original GBPM implementation in \cite{leuzha0801}, the surface differential operator is approximated simply by the polynomial least squares fitting. Although the local approximation accuracy is relatively easy to achieve in local polynomial reconstruction, in some applications, it is more important to construct a discretization on the sampling points so that the resulting linear system can be efficiently inverted stably. For example, when applied to the Laplace-Beltrami operator, one can obtain an M-matrix system using a constrained quadratic programming optimization technique to enforce consistency and diagonal dominance after the discretization \cite{liazha13}. In a recent work \cite{wanleuzha17}, we have developed a modified virtual grid difference (MVGD) discretization. It has led to a more diagonal dominant and better conditioned linear system, which can be efficiently and robustly solved by any existing fast solver for elliptic PDEs. In this work, we propose to follow the CLS approach for approximating local derivatives as developed in \cite{yin17,yinchuleu24}. Numerically, the method can give a better-conditioned discretization of the Laplace-Beltramni operator.

The framework that we are proposing integrates a large class of methods. This allows us to link various current approaches to the same problem. For example, we will also demonstrate that the so-called local tangential lifting (LTL) method \cite{chechiwu15} and a recently developed CLS method with RBF \cite{yin17,yinchuleu24} can be easily incorporated with our framework. 

This paper is organized as follows. We first summarize the GBPM as introduced in \cite{leuzha0801,leuzha0802,leulowzha11} in Section \ref{Sec:GBPM}. We will then introduce a new modification step in the local reconstruction at the end of the Section. In Section \ref{NumMethod}, we introduce our proposed framework for various local reconstructions and explain how the method couples with the original GBPM. In Section \ref{NumericalEx}, we show various numerical experiments to demonstrate the proposed approach's effectiveness and accuracy. Finally, we summarize our proposed method and discuss some possible future works.


\section{Grid Based Particle Method (GBPM)}
\label{Sec:GBPM}

\subsection{The original formulation}

\begin{figure}[!ht]
\begin{center}
\includegraphics[width=0.95\textwidth]{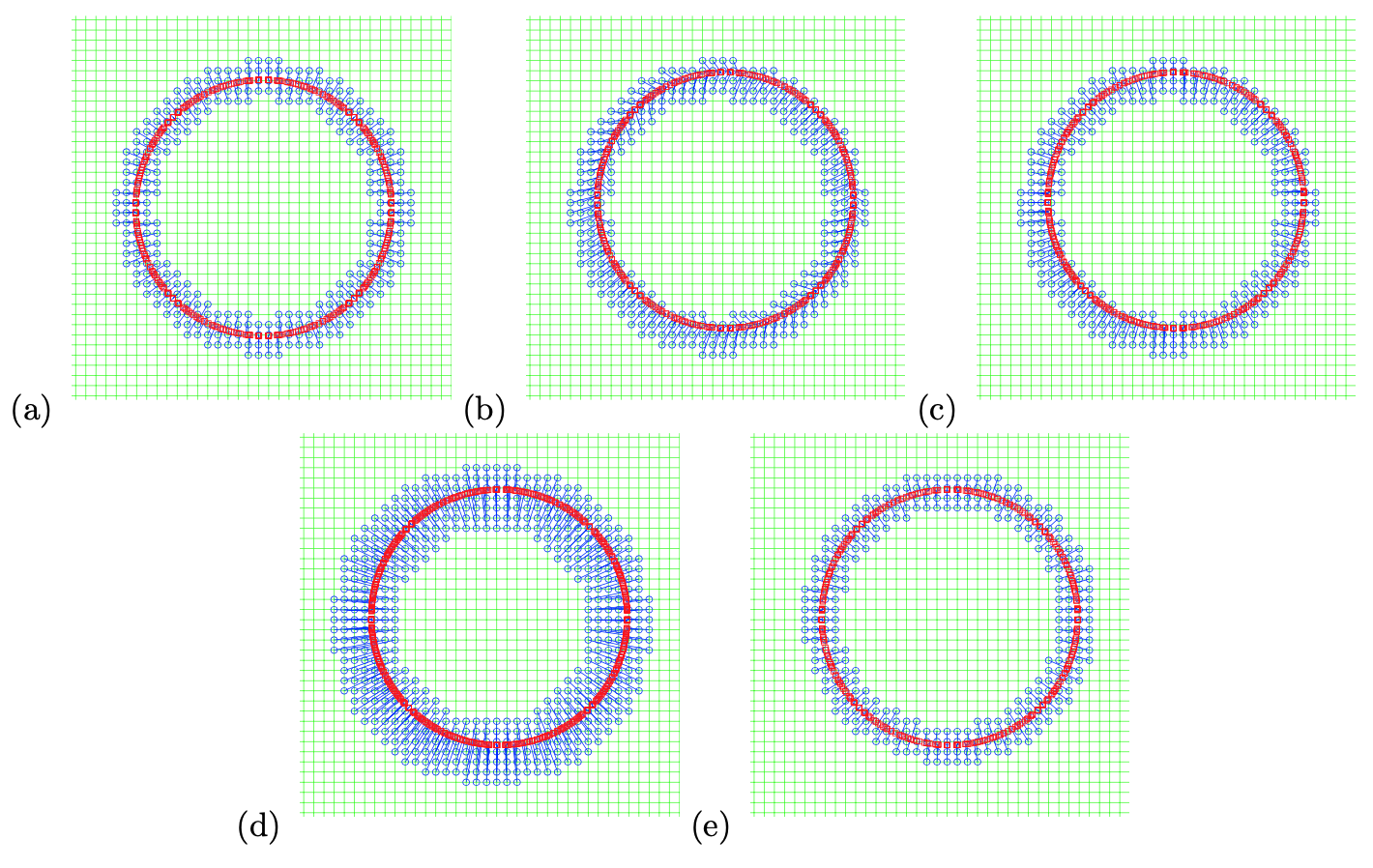}
\end{center}
\caption{Grid Based Particle Method. (a) Initialization, (b) after motion, (c) after re-sampling, (d) after activating new grid points with their footpoints, (e) after inactivating grid points with their footpoints. Blue circles denote active grid points, and red squares denote their corresponding footpoints on the interface. In practice, we do not use such a thick computational tube. We pick it here to show the effect of each step in our algorithm clearly.}
\label{Fig:Algor}
\end{figure}

This section briefly summarizes the Grid Based Particle Method (GBPM) proposed and developed in a series of studies \cite{leuzha0801,leuzha0802,leulowzha11}. We refer the interested readers to the reference and thereafter. The overall algorithm of the GBPM is given as follows, and a graphical illustration is given in Figure \ref{Fig:Algor}.

\vspace{0.5cm} \noindent {\sf Algorithm:
\begin{enumerate}
\item Initialization [Figure \ref{Fig:Algor} (a)]. Collect all grid
points in a small neighborhood (computational tube)
of the interface. From each of these
grid points, compute the closest point on the interface. We call
these grid points \textit{active} and their corresponding particles
on the interface \textit{footpoint}. The interface is represented by
these meshless footpoints (particles).
\item Motion [Figure \ref{Fig:Algor} (b)]. Move all footpoints
according to a given motion law.
\item Re-Sampling [Figure \ref{Fig:Algor} (c)]. For each active
grid point, re-compute the closest point to the interface
reconstructed locally by those particles after the motion in step 2.
\item Updating the computational tube [Figure \ref{Fig:Algor} (d-e)].
Activate any grid point with an active neighbor and
find their corresponding footpoints. Then, inactivate grid points
that are far away from the interface.
\item Adaptation (Optional). Locally refine the underlying grid cell if necessary.
\item Iteration. Repeat steps 2-5 until the final computation time.
\end{enumerate}
}

In the GBPM, we represent the interface by meshless particles associated with an underlying Eulerian mesh. Each sampling particle on the interface is chosen to be the closest point from each underlying grid point in a small neighborhood of the interface. This one-to-one correspondence gives each particle an
Eulerian reference during the evolution. 

In the first step, we define an initial computational tube for active grid points and use their corresponding closest points as the sampling particles for the interface. A grid point is called \textit{active} if its distance to the interface is smaller than a given \textit{tube radius}, and we label the set containing all active grids $\Gamma$. We associate the corresponding closest point on the interface to each active grid point. This particle is called the \textit{footpoint} associated with this active grid point. Furthermore, we can also compute and store certain Lagrangian information of the interface at the footpoints, including normal, curvature, and parametrization, which will be helpful in various applications.

This representation is illustrated in  {Figure} \ref{Fig:Algor}(a) using a circular manifold as an example. We plot the underlying mesh in a solid line, all active grids using small circles, and their associated footpoints on the manifold using squares. To each grid point near the interface (blue circles), we associate a footpoint on the interface (red squares). The relationship of each of these pairs is shown by a solid line link. To track the motion of the interface, we move all the sampling particles according to a given motion law. This motion law can be very general. Suppose the interface is moved under an external velocity field. Since we have a collection of particles on the interface, we move these points just like all other particle-based methods, which is simple and computationally efficient. We can solve a set of ordinary differential equations using a high-order scheme, which gives a very accurate interface location. For more complicated motions, the velocity may depend on the geometry of the interface. This can be done through local interface reconstruction.

It should be noted that a footpoint after motion may not be the closest point on the interface to its associated active grid point anymore. For example,  {Figure} \ref{Fig:Algor} (b) shows the location of all particles on the interface after the constant motion $\textbf{u}=(1,1)^T$ with a small time step. As we can see, these particles on the interface are no longer the closest point from these active grid points to the interface. More importantly, the motion may cause those original footpoints to become unevenly distributed along the interface. This may introduce stiffness when particles are  {moving} toward each other and a large error when particles are  {moving} apart. To maintain a quasi-uniform distribution of particles, we need to resample the interface by recomputing the footpoints and updating the set of active grid points ($\Gamma$) during the evolution. During this resampling process, we locally reconstruct the interface, which involves communications among different particles on the interface. This local reconstruction also provides geometric and Lagrangian information at the recomputed footpoints on the interface. The key step in the method is a least squares approximation of the interface using polynomials at each particle in a local coordinate system, $\{ (\textbf{n}')^{\perp}, \textbf{n}' \}$, with the associated particle as the origin. Using this local reconstruction, we find the closest point from this active grid point to the local approximation of the interface. This gives the new footpoint location. Further, we also compute and update any necessary geometric and Lagrangian information, such as the normal vector, curvature, and possibly an updated interface parametrization at this new footpoint. 

\subsection{A new weighting strategy in the local reconstruction}
\label{MLocalRec}

We approximate the surface locally using a polynomial from each grid point $\mathbf{p}$. For example, in the three dimensions, we could use a second-order polynomial given by

$$
h(x,y) = a_{0,0} + a_{1,0}x+ a_{0,1}y + a_{2,0} x^2 + a_{1,1}xy + a_{0,2} y^2 \, ,
$$
where $(x,y,h(x,y))$ are coordinates in the local coordinate system, $\{ (\textbf{n}')^{\perp}, \textbf{n}' \}$. The original idea in \cite{leuzha0801} fits the polynomial using a collection of $(n+1)$-footpoints, $(x_i,y_i,h(x_i,y_i))$ with $i = 0,1, \cdots, n$, in the neighborhood of $\mathbf{p}$. With these points, we determine the least squares solution of the linear system:

$$
\begin{pmatrix}
1 & x_0 & y_0 & x_0^2 & x_0 y_0 & y_0^2\\
1 & x_1 & y_1 & x_1^2 & x_1 y_1 & y_1^2\\
1 & x_2 & y_2 & x_2^2 & x_2 y_2 & y_2^2\\
\vdots & \vdots & \vdots & \vdots & \vdots & \vdots\\
1 & x_n & y_n & x_n^2 & x_n y_n & y_n^2
\end{pmatrix}\begin{pmatrix}
a_{0,0}\\
a_{1,0}\\
a_{0,1}\\
a_{2,0}\\
a_{1,1}\\
a_{0,2}
\end{pmatrix} = \begin{pmatrix}
h(x_0,y_0)\\
h(x_1,y_1)\\
h(x_2,y_2)\\
\vdots\\
h(x_n,y_n)
\end{pmatrix} \, .
$$
This set of footpoints in the local reconstruction must be chosen carefully. As discussed in the original GBPM formulation, this set of points has to satisfy the following conditions. The first one is that these points have to be distinct. In the previous formulation, we have conditioned that these footpoints must be separated by at least a certain distance on the local tangent plane. This constraint avoids putting an extensive amount of weight on these sampling points. The second constraint is that we would like to sample the same segment on the interface. When the curvature of the interface is too large compared to the underlying mesh size, one might collect footpoints from different interface segments, which will result in a wrong local reconstruction of the surface. In the original work, we proposed checking if the associated normal vectors from these footpoints were pointing in a similar direction. In particular, for each of these footpoints, $\x_i$ has an associated vector $\n'_i$ corresponding to the surface normal before the motion. We have checked if $\n'_i\cdot \n'_j > \cos\theta_{\max}$ for some angle $\theta_{\max}$. 

The above strategies are in some sense imposing some \textit{hard constraints} in the selection phase. We either use the neighboring sampling point or eliminate it when choosing the subset of footpoints in the local reconstruction step. This paper proposes relaxing this by imposing some \textit{soft constraints}. The first condition on the well-separation can be partly handled by the moving least squares (MLS) method \cite{lansal81}, which introduces a weighting to the footpoints based on the distance between points on the tangent plane. We will not concentrate on this part. To avoid information from different segments on the surface, instead of simply getting rid of those points pointing at an angle larger than $\theta_{\max}$, we propose to add a weighting in the least squares further to emphasize the contribution from footpoints points with similar normal directions. Mathematically, we let $\x_0 = (x_0,y_0,h(x_0,y_0))$ be the closest footpoint to the grid point $\mathbf{p}$ and denote those normal vectors at the footpoints before motion by $\n'_i, i = 0,1, \cdots, n$. We now introduce a weight $c_i$ to each collected footpoints given by

\begin{equation}
c_i=\left\{
\begin{array}{cc}
\n'_0\cdot \n'_i & \mbox{ if $\n'_0\cdot \n'_i > \cos\theta_{\max}$,} \\
0 & \mbox{ otherwise.}
\end{array}
\right. \label{Eqn:Weighting}
\end{equation}
Once we have this set of weights, the $i$-th equation in the least squares system will be modified to

$$
c_i \begin{pmatrix}
1 & x_i & y_i & x_i^2 & x_i y_i & y_i^2
\end{pmatrix}\begin{pmatrix}
a_{0,0}\\
a_{1,0}\\
a_{0,1}\\
a_{2,0}\\
a_{1,1}\\
a_{0,2}
\end{pmatrix} = c_i \, h(x_i,y_i)  \, .
$$



\section{Our proposed framework}
\label{NumMethod}

This section discusses our proposed framework for solving PDEs on an evolving surface. We will first discuss the mathematical background of the problem and discuss all the necessary quantities we have to compute. Then, we introduce our CLS method approach for local approximation.

\subsection{Background}

Let $\Gamma$ be a regular surface in $\mathbb{R}^2$ parametrized by $\x$. Mathematically, the surface gradient of a smooth function $\xi$ is given by

$$
\nabla_{\Gamma} \xi = \frac{\xi_u G-\xi_v F}{EG-F^2}\x_u + \frac{\xi_v E-\xi_u F}{EG-F^2}\x_v,
$$
where $E, F$ and $G$ are the coefficients of the first fundamental form of the regular surface $\Gamma$, the quantities $\xi_u = \frac{\partial}{\partial u} \xi(\x(u,v))$ and $\xi_v = \frac{\partial}{\partial v} \xi(\x(u,v))$. For any local vector field $X = A\x_u+B\x_v$, we have the surface divergence defined as

$$
\nabla _{\Gamma} \cdot X = \frac{1}{\sqrt{EG-F^2}}\left[ \frac{\partial}{\partial u }\left(A\sqrt{EG-F^2}\right)+ \frac{\partial}{\partial v}\left(B\sqrt{EG-F^2}\right)\right] \, .
$$
The surface Laplacian (Laplace-Beltrami) operator $\Delta_{\Gamma}\xi$ is defined as $\Delta_{\Gamma}\xi = \nabla _{\Gamma} \cdot(\nabla_{\Gamma}\xi)$. In the local representation, it can be expressed as
\begin{eqnarray*}
\Delta_{\Gamma} \xi &=& \frac{1}{\sqrt{EG-F^2}}\left[ \frac{\partial}{\partial u }\left(\frac{G}{\sqrt{EG-F^2}} \xi_u\right)-\frac{\partial}{\partial u }\left(\frac{F}{\sqrt{EG-F^2}} \xi_v\right) \nonumber \right.\\
&& + \left.\frac{\partial}{\partial v }\left(\frac{E}{\sqrt{EG-F^2}} \xi_v\right)-\frac{\partial}{\partial v }\left(\frac{F}{\sqrt{EG-F^2}} \xi_u\right)\right] \, .
\end{eqnarray*}
Therefore, to design a numerical method for computing any differential operator on the surface, we have to approximate the first fundamental form of the surface and also provide a way to approximate the derivatives in the local coordinate systems.

\begin{figure}[!ht]
\centering
\includegraphics[width=0.5\textwidth]{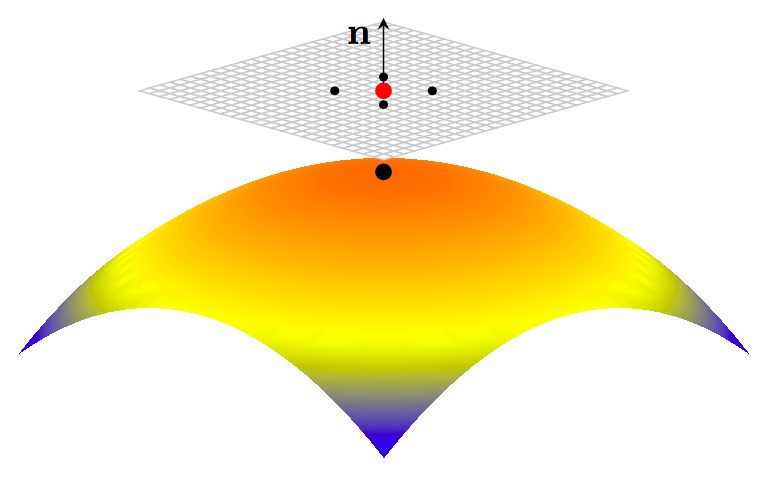}
\caption{The local reconstruction step involves a construction of the tangent plane using a normal vector $\n$ estimated from the data and a projection step of sampled points collected from a local neighborhood.}
\label{Fig:LocalReconstruction}
\end{figure}

\subsection{The CLS-GSP reconstruction}

The first issue is to design a good local representation of surface $\Gamma$. One natural representation is the local tangential projection as shown in Figure \ref{Fig:LocalReconstruction}. We assume that a set of sampling points represents the interface. These sampling points can be obtained from various approaches. For example, it can be the GBPM representation we summarized above, a typical point cloud dataset, or an interface represented by triangulation. The representation itself does not play a crucial role in the following discussions. Assume further that for any individual sampling point $\mathbf{y}^*$, we can collect a set of $M$-neighboring sampling points denoted by $\{\mathbf{y}_j: j = 1, \cdots, M\}$.  It is certainly easier to have the interface represented by a triangulation since we have the connectivities of all vertices (the sampling points) on the interface. For the GBPM representation, we can collect these neighboring points through the underlying mesh. Some KNN algorithms can also be applied to point cloud datasets. 

For the sampling point $\mathbf{y}^*$, we let $\n$ be the corresponding (approximated) normal vector of the surface. With $\n$ and its tangent plane, we introduce a new local coordinate system $(x,y,z)$ with the $z$-direction representing the component in the $\n$-direction and $\mathbf{y}^*$ becomes $(\x_0,0)$ in the new coordinate system. In this local coordinate system, a function can locally approximate the surface, and we denote it as $h(x_j,y_j) = z_j$. Once we have this local approximation $h$, we can approximate the derivatives $h_x, h_y, h_{xx}, h_{xy}$ and $h_{yy}$ for the first fundamental form in the surface differential operators. 

Numerically, this leads to one of the following problems. The first one is an interpolation problem: Given $(x_j,y_j,z_j)$ for $j=1,\cdots,n$, determine an interpolant $h(x,y)$ which passes through all data points such that $(x_j,y_j,h(x_j,y_j))$. One idea is to use the polynomial interpolation. For example, when the neighbors are carefully chosen such that $n=6$, we can use a second-order polynomial to interpolate these data points. This approach better estimates the local error if the interpolation uses the polynomial basis. However, we generally have no control over the distribution of $(x_j,y_j)$ from the projection of the sampling points onto the local tangent plane. In particular, some of these projection points could be very close to each other or could collide on the tangent plane. Designing a good subset of points might be challenging, which could lead to a \textit{good} interpolation problem. Therefore, we do not consider this class of approach further.

The second is the least squares approximation problem: Given $(x_j,y_j,z_j)$ for $j=1,\cdots,n$, determine a function $h(x,y)$ that best explains the data points. Mathematically, we obtain the set of parameters $\alpha_i$ to minimize the mismatch between the data points and $\sum_{\psi_i \in \mathcal{B}} \alpha_i \psi_i(\x)$ where $ \mathcal{B}$ is the space for the basis, i.e. we solve

$$
\min_{\alpha_i} \sum_{j=0}^M \left[ h(\x_j) - \sum_{\psi_i \in \mathcal{B}} \alpha_i \psi_i(\x_j) \right]^2 \, .
$$
The choice of the basis is clearly not unique. For example, the original GBPM developed in \cite{leuzha0801,leuzha0802,leulowzha11} chose the standard polynomial basis for $\mathcal{B}$. Of course, one could also replace the standard least squares approximation with a weighted least squares approach so that those points far away from $\x_0$ do not contribute equally like the point $\x_0$ itself.

In a recent work \cite{yin17,yinchuleu24}, we have developed a general strategy to impose a hard constraint in the least squares reconstruction by enforcing the least squares reconstruction passing through $(\x_0, h(\x_0))$ exactly. We call this a CLS approach. Under this constraint, the least squares method can be modified to fit the function value $h(\x) - h(\x_0)$ with functions in the form of $\sum_{\psi_i \in \mathcal{B}} \beta_i [\psi_i(\x)- \psi_i(\x_0)]$. This only requires a simple modification of the standard least squares approximation method. We will give out more information below. Let $\mathcal{L}$ be some linear differential operator. We can also follow a similar idea to obtain an approximation $\mathcal{L}(h(\x_0))$ based on the CLS, i.e.

$$
\mathcal{L}(h(\x_0)) = \sum_{j=1}^M w_j (h(\x_j)-h(\x_0)) = \sum_{j=0}^M w_j h(\x_j) \, ,
$$
with the coefficient $w_0 = -\sum_{j=1}^M w_j$.

The choice of the basis is not unique. In this paper, we present two different sets of bases and give the technical implementation details. The first is the simple polynomial basis, which leads to an approximation result similar to the original GBPM method for PDEs but with an extra hard constraint at the target location during reconstruction. The second basis follows our previous development in \cite{yin17,yinchuleu24}, which uses the radial basis functions. A different basis will give slightly different results. We are not trying to conclude which of the following two basis sets will yield better numerical results in all applications. Still, we will only show some examples to explain better the framework that we are proposing.

\subsubsection{Standard polynomial basis}

The first set of basis in consideration is the standard polynomial basis. The weights $w_j$ in the numerical approximation $\mathcal{L}(h(\x_0)) = \sum_{j=0}^M w_i h(\x_j)$ satisfy

$$
\left( w_1, w_2, \cdots, w_M  \right) =\left( \left\lbrace \mathcal{L} P_j(\x)\Big|_{\x = \x_0} \right\rbrace_{j=1}^k \right) C^+ \, ,
$$
where $k$ depends on the order of those polynomials in the basis $\mathcal{B}$, the matrix $C$ is the corresponding polynomial basis coefficients, and $C^+$ is the pseudo-inverse of the matrix $C$, which can be easily computed using the SVD or the QR algorithm. We only use the second-order least squares polynomial fitting at $\x_0 = \mathbf{0}$. A similar formulation can be easily extended to any higher-order polynomial fitting. Set
$$
G_1 = \begin{pmatrix}
\left\lbrace \partial_x P_j(\x)\Big|_{\x = \mathbf{0}} \right\rbrace_{j=1}^k\\
\left\lbrace \partial_y P_j(\x)\Big|_{\x = \mathbf{0}} \right\rbrace_{j=1}^k\\
\left\lbrace \partial_{xx} P_j(\x)\Big|_{\x = \mathbf{0}} \right\rbrace_{j=1}^k\\
\left\lbrace \partial_{xy} P_j(\x)\Big|_{\x = \mathbf{0}} \right\rbrace_{j=1}^k\\
\left\lbrace \partial_{yy} P_j(\x)\Big|_{\x = \mathbf{0}} \right\rbrace_{j=1}^k
\end{pmatrix} = \begin{pmatrix}
1 & 0 & 0 & 0 & 0\\
0 & 1 & 0 & 0 & 0\\
0 & 0 & 2 & 0 & 0\\
0 & 0 & 0 & 1 & 0\\
0 & 0 & 0 & 0 & 2\\
\end{pmatrix},
$$
where $P_j(\x) \in \{x, y, x^2, xy, y^2\}$. Then we have

$$
\begin{pmatrix}
w_1^1 & w_2^1 & \cdots & w_M^1 \\
w_1^2 & w_2^2 & \cdots & w_M^2 \\
w_1^3 & w_2^3 & \cdots & w_M^3 \\
w_1^4 & w_2^4 & \cdots & w_M^4 \\
w_1^5 & w_2^5 & \cdots & w_M^5 
\end{pmatrix} =G_1 \begin{pmatrix}
x_1 & y_1 & x_1^2 & x_1 y_1 & y_1^2\\
x_2 & y_2 & x_2^2 & x_2 y_2 & y_2^2\\
\vdots & \vdots & \vdots & \vdots & \vdots\\
x_M & y_M & x_M^2 & x_M y_M & y_M^2
\end{pmatrix}^+ .
$$

This numerical scheme is the so-called local tangential lifting (LTL) developed in \cite{chechiwu15}. Some nice properties can be proved under the standard polynomial basis. We do not give the full numerical details of this approach but refer the interested readers to the reference and after that. 

\subsubsection{Radial basis functions}

An interesting set of basis is the radial basis functions (RBF). Since the projection of the sampling points from the local neighborhood could collide on the tangent plane, this could sometimes create numerical instability in the standard RBF interpolation. In this paper, we follow \cite{yin17,yinchuleu24} and apply the CLS-GSP method to approximate the local surface and any function defined on the evolving manifold.

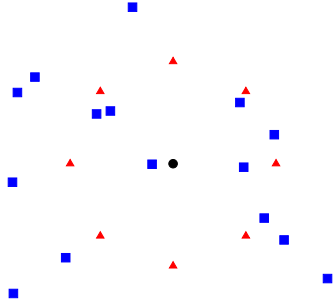
\begin{figure}[!ht]
\centering
\begin{tikzpicture}[scale=0.80]
	\begin{axis}[%
	hide axis, xtick=\empty, ytick=\empty,axis equal,
	scatter/classes={%
		a={mark=square*,blue},%
		b={mark=triangle*,red},%
		c={mark=*,draw=black}}]
	\addplot[scatter,only marks,%
		scatter src=explicit symbolic]%
	table[meta=label] {
x     y      label
0 0 c
0.1 0 b
-0.1 0 b
0 -0.1 b
0 0.1 b
0.0707 0.0707 b
-0.0707 0.0707 b
0.0707 -0.0707 b
-0.0707 -0.0707 b
0.0686  -0.0034 a
-0.1560  -0.0181 a
-0.0743   0.0487 a
-0.1512   0.0697 a
-0.1342   0.0848 a
0.1078  -0.0746 a
0.0649   0.0599 a
-0.0609   0.0516 a
0.1500  -0.1124 a
-0.1551  -0.1270 a
-0.0204  -0.0005 a
-0.0394   0.1532 a
0.0885  -0.0532 a
0.0983   0.0284 a
-0.1043  -0.0920 a
	};
	\end{axis}
\end{tikzpicture}
\caption{The local samplings of the function. The black dot denotes the target point $\x_0$. The blue squares are the $n$-nearest scattered data points projection onto the local tangent plane at $\x_0$. The red triangles are the ghost sample points surrounding the target point $\x_0$ on the tangent plane.}
\label{Fig:GhostPoints}
\end{figure}

We propose to convert the interpolation problem in the standard RBF method to a regression problem on a set of prescribed sample locations. A simple two-dimensional demonstration is shown in Figure \ref{Fig:GhostPoints}. Because we do not require the interpolating values at all sampling points to match the function values, we relax the constraint by minimizing the mismatch at a different set of $d$ sampling locations in the least squares sense. We call this set of new sampling locations the \textit{ghost sample points} (as shown using the red triangles) and denote them by $\{\hat{\x}_i\}$. The choice of these ghost sample points is not unique. Here, for any target point $\x_0$, we propose to pick $d$ ghost sample points uniformly surrounding it. In this paper, we collect 16 local neighbors for each sampling point in various local reconstructions (i.e., $n=16$). In Figure \ref{Fig:GhostPoints}, the projection of these sampling points onto the local tangent plane is plotted using blue squares and the black circle (representing the target location $\x_0$ itself). Then we apply 8 ghost sample points (i.e. $d=8$) distributed uniformly at a distance $r= \frac{1}{2}\max_i \|\x_0-\x_i\|$ from each individual sampling point $\x_0$ (plotted by the red triangles).

Following the idea of the CLS, we enforce the least squares reconstruction passing through the function value $(\x_0, h(\x_0))$ exactly. Mathematically, this constraint implies 

$$
h(\mathbf{0}) = \sum_{i=1}^d \lambda_i \psi(\| \hat{\x}_i \|) +\sum_{j=0}^{k} \mu_j P_j(\mathbf{0}) = \sum_{i=1}^d \lambda_i \psi(\| \hat{\x}_i \|) +\mu_0 \, ,
$$
where $P_j$ is the polynomial basis. If we define 

$$
v(\x) = h(\x)-h(\mathbf{0}) = \sum_{i=1}^d \lambda_i \left(\psi(\| \x-\hat{\x}_i \|)-\psi(\| \hat{\x}_i \|) \right) +\sum_{j=1}^{k} \mu_j P_j(\x) \, ,
$$
we find that this function $v(\mathbf{0})\equiv 0$ for any choice of $\lambda_i$ and $\mu_j$. In other words, any function in the form of

$$
\sum_{i=1}^d \lambda_i \left(\psi(\| \x-\hat{\x}_i \|)-\psi(\| \hat{\x}_i \|) \right) +\sum_{j=1}^{k} \mu_j P_j(\x) + h(\mathbf{0})
$$
passes through $(\x_0, h(\x_0))$ exactly for any $\lambda_i$ and $\mu_j$. Then, we follow an idea similar to the LS-RBF to determine the coefficients $\lambda_i$ and $\mu_j$ in the least squares sense. We redefine a new kernel function by

$$
\varphi(\|\x - \hat{\x}_i \|) =  \psi(\| \x-\hat{\x}_i \|)-\psi(\| \hat{\x}_i \|) \, .
$$
Then, the coefficients can be determined by obtaining the least squares solution to the following system of linear equations
\begin{eqnarray}
\begin{pmatrix}
D & \begin{matrix}
x_1 & y_1 & \cdots\\
\vdots & \vdots &\\
x_M & y_M & \cdots
\end{matrix} \\
\begin{matrix}
\hat{x}_1 & \cdots & \hat{x}_d\\
\hat{y}_1 & \cdots & \hat{y}_d\\
\vdots & & \vdots
\end{matrix}&0
\end{pmatrix}  \begin{pmatrix}
\lambda_1\\
\vdots\\
\lambda_d\\
\mu_1\\
\vdots\\
\mu_k
\end{pmatrix} = \begin{pmatrix}
v(\x_1)\\
\vdots\\
v(\x_M)\\
0\\
\vdots\\
0
\end{pmatrix},
\label{Eqn.RegularizedRBFFitting}
\end{eqnarray}
where $D_{p, j} = \varphi (\|x_p - \hat{x}_i\|)$ for $p = 1, \cdots, M$ and $i = 1, \cdots , d$. 

Further, to obtain an approximation to the differential operator, we first let $\hat{D}$ be the coefficient matrix on the left hand side of the equation (\ref{Eqn.RegularizedRBFFitting}). Then we have
\begin{eqnarray}
\left( \hat{w}_1,\cdots, \hat{w}_M, \hat{w}_{M+1},\cdots ,\hat{w}_{M+k}  \right) =\left( \left\lbrace \mathcal{L}\varphi(\| \x-\hat{\x}_i \|)\Big|_{\x=\x_0} \right\rbrace_{i=1}^d , \left\lbrace \mathcal{L} P_j(\x)\Big|_{\x = \x_0} \right\rbrace_{j=1}^k \right) \hat{D}^+ \, ,
\label{Eqn.WeightRRBFFitting}
\end{eqnarray}
where $\hat{D}^+$ denotes the pseudo-inverse of the matrix $\hat{D}$. Note that

$$
\mathcal{L}\varphi(\| \x-\hat{\x}_i \|)\Big|_{\x=\x_0} = \mathcal{L}\phi(\| \x-\hat{\x}_i \|)\Big|_{\x=\x_0}
$$ 
for any differential operator. With the weights obtained, we have

$$
\mathcal{L}h(\x)\Big|_{\x=\x_0} = \sum_{p=1}^M \hat{w}_p v(\x_p) = \sum_{p=1}^M \hat{w}_p h(\x_p) - \left( \sum_{p=1}^M \hat{w}_p\right) h(\x_0) \, .
$$
In particular, if we use the Gaussian kernel given by $\phi(r) = \exp[{-(\varepsilon r)^2}]$ and set

$$
G_2 = \begin{pmatrix}
\left\lbrace \partial_{x}\varphi(\| \x-\hat{\x}_i \|)\Big|_{\x=\mathbf{0}} \right\rbrace_{i=1}^d\\
\left\lbrace \partial_{y}\varphi(\| \x-\hat{\x}_i \|)\Big|_{\x=\mathbf{0}} \right\rbrace_{i=1}^d\\
\left\lbrace \partial_{xx}\varphi(\| \x-\hat{\x}_i \|)\Big|_{\x=\mathbf{0}} \right\rbrace_{i=1}^d\\
\left\lbrace \partial_{xy}\varphi(\| \x-\hat{\x}_i \|)\Big|_{\x=\mathbf{0}} \right\rbrace_{i=1}^d\\
\left\lbrace \partial_{yy}\varphi(\| \x-\hat{\x}_i \|)\Big|_{\x=\mathbf{0}} \right\rbrace_{i=1}^d
\end{pmatrix} = \begin{pmatrix}
\left\lbrace 2 \hat{x}_i \varepsilon^2 \exp[{-(\varepsilon \| \hat{\x}_i\|)^2}] \right\rbrace_{i=1}^d\\
\left\lbrace 2 \hat{y}_i \varepsilon^2 \exp[{-(\varepsilon \| \hat{\x}_i\|)^2}] \right\rbrace_{i=1}^d\\
\left\lbrace (-2\varepsilon^2 + 4\hat{x}_i^2 \varepsilon^4) \exp[{-(\varepsilon \| \hat{\x}_i\|)^2}] \right\rbrace_{i=1}^d\\
\left\lbrace 4 \hat{x}_i \hat{y}_i \varepsilon^4 \exp[{-(\varepsilon \| \hat{\x}_i\|)^2}] \right\rbrace_{i=1}^d\\
\left\lbrace (-2\varepsilon^2 + 4\hat{y}_i^2 \varepsilon^4) \exp[{-(\varepsilon \| \hat{\x}_i\|)^2}] \right\rbrace_{i=1}^d\\
\end{pmatrix} \, ,
$$
we can approximate the coefficients by

$$
\begin{pmatrix}
\hat{w}_1^1 & \hat{w}_2^1 & \cdots & \hat{w}_M^1 & \hat{w}_{M+1}^1 & \cdots \hat{w}_{M+k}^1\\
\hat{w}_1^2 & \hat{w}_2^2 & \cdots & \hat{w}_M^2 & \hat{w}_{M+1}^2 & \cdots \hat{w}_{M+k}^2\\
\hat{w}_1^3 & \hat{w}_2^3 & \cdots & \hat{w}_M^3 & \hat{w}_{M+1}^3 & \cdots \hat{w}_{M+k}^3\\
\hat{w}_1^4 & \hat{w}_2^4 & \cdots & \hat{w}_M^4 & \hat{w}_{M+1}^4 & \cdots \hat{w}_{M+k}^4\\
\hat{w}_1^5 & \hat{w}_2^5 & \cdots & \hat{w}_M^5 & \hat{w}_{M+1}^5 & \cdots \hat{w}_{M+k}^5
\end{pmatrix} =\begin{pmatrix}
G_2 & G_1\\
\end{pmatrix} \hat{D}^+ \, .
$$
Omitting the last $k$ columns of coefficients, we have

$$
\begin{pmatrix}
h_x(\x_0)\\
h_y(\x_0)\\
h_{xx}(\x_0)\\
h_{xy}(\x_0)\\
h_{yy}(\x_0)
\end{pmatrix} = \begin{pmatrix}
\hat{w}_1^1 & \hat{w}_2^1 & \cdots & \hat{w}_M^1\\
\hat{w}_1^2 & \hat{w}_2^2 & \cdots & \hat{w}_M^2\\
\hat{w}_1^3 & \hat{w}_2^3 & \cdots & \hat{w}_M^3\\
\hat{w}_1^4 & \hat{w}_2^4 & \cdots & \hat{w}_M^4\\
\hat{w}_1^5 & \hat{w}_2^5 & \cdots & \hat{w}_M^5 
\end{pmatrix} \begin{pmatrix}
h(\x_1)-h(\x_0)\\
h(\x_2)-h(\x_0)\\
\vdots\\
h(\x_M)-h(\x_0)
\end{pmatrix} \, .
$$
Finally, the partial derivatives ${\partial}/{\partial x}$, ${\partial}/{\partial y}$, ${\partial^2}/{\partial x^2}$, ${\partial^2}/{\partial y^2}$ and ${\partial^2}/{\partial x \partial y}$ can be easily constructed using the same set of coefficients. Further, with both $h_x$ and $h_y$ approximated, we can determine the coefficients in the first fundamental form using
$$
E = 1+ h_x^2 \mbox{ , } F = h_x h_y \mbox{ and } G = 1+h_y^2 \, .
$$

Note that the coefficients $\hat{w}_j$ depend only on the distribution of the projected points on the local tangent place but \textit{not} on the function values of $h$. We can apply the same set of coefficients $\hat{w}_j$ to compute $\xi_x, \xi_y, \xi_{xx}, \xi_{xy}$ and $\xi_{yy}$ at $\x_0$. In particular, we have

$$
\nabla_{\Gamma} \xi = \frac{\left(\xi_x(1+h_x^2)-\xi_y h_x h_y, \xi_y(1+h_y^2)-\xi_x h_x h_y, \xi_x \xi_y \right)}{1+h_x^2+h_y^2} \, .
$$

\begin{rem}
Indeed, to get a more accurate approximation of the surface differential operator, we need to collect more neighboring points for each sample point. If the sampling density of the interface cannot grow with accuracy, we have to go further away to collect enough neighbors in the least squares approximation. However, this will lower the accuracy of the local approximation of the interface. Therefore, we need to assume that the number of sampling points has to be large enough, and the sampling density must be uniform.
\end{rem}

\begin{rem}
In the discussions above, we only require the fitting function to pass through the target point exactly. Except for the target point, the remaining samples are treated equally with the same weighting in the least squares fitting. However, as discussed in Section \ref{MLocalRec}, we introduce an extra weighting to emphasize the contributions from sample points based on the angle of its normal direction compared to the target one.
\end{rem}

\subsection{The overall algorithm for solving PDEs on evolving surfaces}
\label{CoupleMethod}

To end this section, we summarize our overall approach to solving PDEs on evolving surfaces. The original GBPM handles both the interface sampling and the interface motion. At the same time, the proposed CLS methods help to approximate the surface, the local reconstructions of the function value $u$ on the interface, and the coefficients in the first fundamental form for the surface differential operator.

\vspace{0.5cm} \noindent {\sf Algorithm:
\begin{enumerate}
\item Initialization.
\item Generate the surface operators.
\begin{enumerate}
\item Determine neighboring grids (inside the computational band) and collect their corresponding footpoints. Denote the number of points by $M$.
\item Construct the local coordinate system at each sampling point.
\item Apply the CLS method to construct the weights $w_j$ for $j=1,\cdots,M$.
\item Generate the surface differential operators.
\end{enumerate}
\item Update the function value by solving the time-dependent surface PDEs.
\item Update the location of foot points according to the motion law.
\item Follow the GBPM local reconstruction, including resampling new footpoints and updating the corresponding function value.
\item Back to Step 2, repeat until the stopping time.
\end{enumerate}
}

\section{Numerical experiments}
\label{NumericalEx}

In this section, we present several numerical tests to demonstrate the effectiveness of the proposed method. We will first solve some PDEs on static surfaces represented by the GBPM to show the performance of the CLS. Then, we add different motion laws to the interface for various applications. The examples in this section will be arranged in the following way. We will first look at the effect on the new weightings in the local reconstruction in Section \ref{Ex:MGBPM} without consideration of any PDE on the evolving surface. Then, we concentrate on solving PDEs on a static surface sampled by the GBPM representation in Section \ref{SubSec:RLS-GBPM-Static}. This section tests the ability and effectiveness of the newly proposed CLS-GSP method. We consider in Section \ref{Ex:Convergence} a simple example for solving an advection-diffusion equation on an evolving surface where the exact solution to the problem can be found. We will then present some simulations of various equations in Section \ref{Ex:CoupledEvolvingTorus} to Section \ref{Ex:CahnHilliardEllipsoid} to demonstrate the effectiveness of the proposed approach.


\begin{figure}[!ht]
\begin{center}
\includegraphics[width=0.95\textwidth]{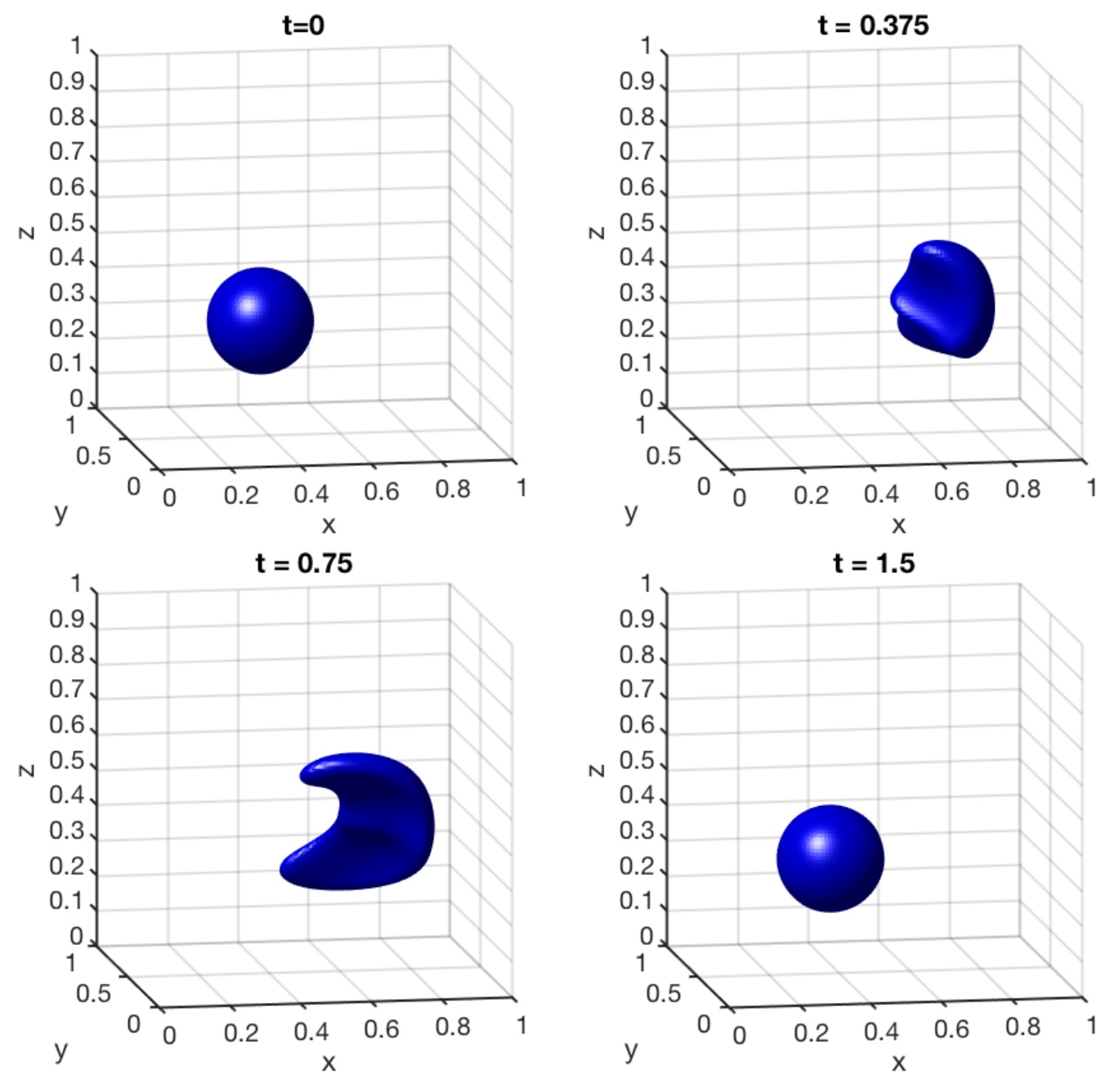}
\end{center}
\caption{(Example \ref{Ex:MGBPM}) Motion under a single vortex with the reversible motion using the resolution $100^3$ at $t = 0, 0.375, 0.75$, and 1.5.}
\label{Fig.MGBPM}
\end{figure}

\begin{figure}[!ht]
\centering
\includegraphics[width=0.95\textwidth]{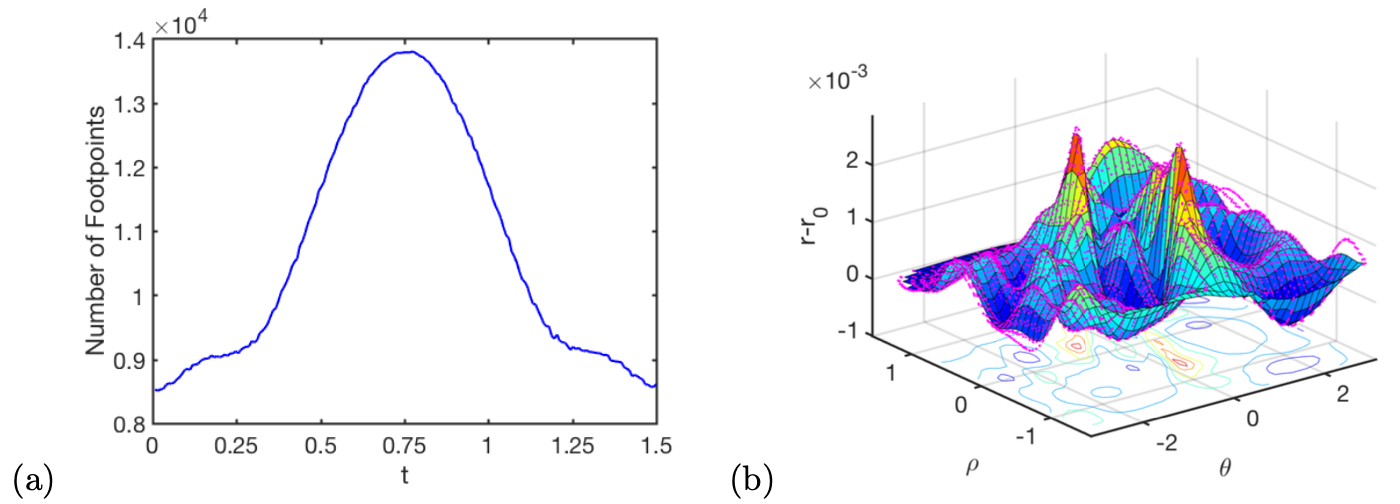} 
\caption{(Example \ref{Ex:MGBPM}) (a) The change in the number of footpoints at different times and (b) the error in the radius of the sphere $r-r_0$ at the final time.}
\label{Fig.MGBPMErr}
\end{figure}

\begin{figure}[!ht]
\begin{center}
\includegraphics[width=0.95\textwidth]{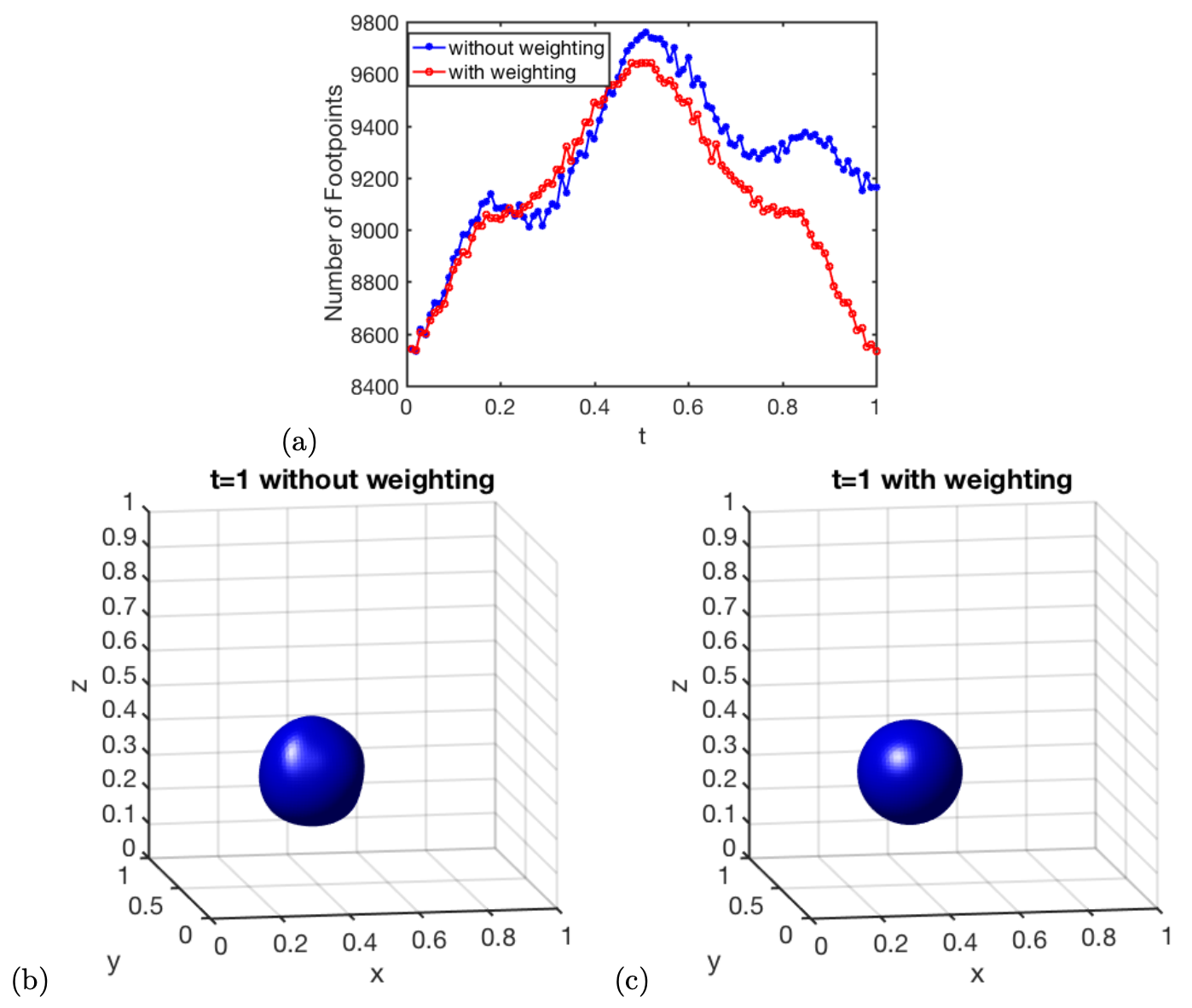}
\end{center}
\caption{(Example \ref{Ex:MGBPM}) (a) Number of footpoints during the evolution. (b) The surface after the final stage of the evolution without weighting in local reconstruction. (c) The surface after the final stage of the evolution with weighting in local reconstruction.}
\label{Fig.MGBPMCompare}
\end{figure}

\subsection{The modified GBPM}
\label{Ex:MGBPM}

Our first example is the motion of a sphere under a single vortex flow \cite{effm02,hiekou05,mingib07,leuzha0801}. This is a challenging test case that can test the stability of our GBPM with the modified local reconstruction. The velocity field is given by
\begin{eqnarray*}
v_1(x,y,z) &=& 2\sin(\pi x)^2 \sin(2\pi y) \sin(2\pi z) \cos\left(\frac{\pi t}{T}\right)\\
v_2(x,y,z) &=& -\sin(2\pi x)\sin(\pi y)^2 \sin(2\pi z) \cos\left(\frac{\pi t}{T}\right)\\
v_3(x,y,z) &=& -\sin(2\pi x) \sin(2\pi y) \sin(\pi z)^2 \cos\left(\frac{\pi t}{T}\right) \, ,
\end{eqnarray*}
for some given final time $T$. The initial shape is a sphere centered at $(0.35,0.35,0.35)$ with radius $r_0 = 0.15$. In Figure \ref{Fig.MGBPM}, we compute the solution at $t = 0, 0.375, 0.75$ and $1.5$ with $T=1.5$ using $\Delta x = 0.01$. The mean distance from the final footpoints to the point $(0.35,0.35,0.35)$ at the final time is $0.1506$ with the standard deviation $5.114\times 10^{-4}$. The largest absolute error in the location is $2.862 \times 10^{-3}$. There is a total of 
$8530$ initially at $t=0$. The number gradually increases to the maximum value $13801$ when $t=0.75$, finally dropping to $8615$ at the final time. Figure \ref{Fig.MGBPMErr} shows the detailed evolution of footpoints and the error in the sphere's radius $r-r_0$ at the final time.

To see the effect of the new weighting (\ref{Eqn:Weighting}) introduced in the local reconstruction, we compare two numerical solutions in Figure \ref{Fig.MGBPMCompare}. In this example, we repeat the similar setup but choose a slightly smaller final time with $T=1$. Figure \ref{Fig.MGBPMCompare} shows the changes in the number of footpoints and the solution at the final time $T=1$. For the original scheme without any weighting in the local reconstruction, we found that the mean distance from the final footpoints to the center at the final time is $0.1557$ with the standard deviation $5.832\times 10^{-3}$. The modified local reconstruction gives a more accurate solution with the corresponding mean distance $0.1502$ with the standard deviation $2.247 \times 10^{-4}$.


\begin{figure}[!ht]
\centering
\includegraphics[width=0.95\textwidth]{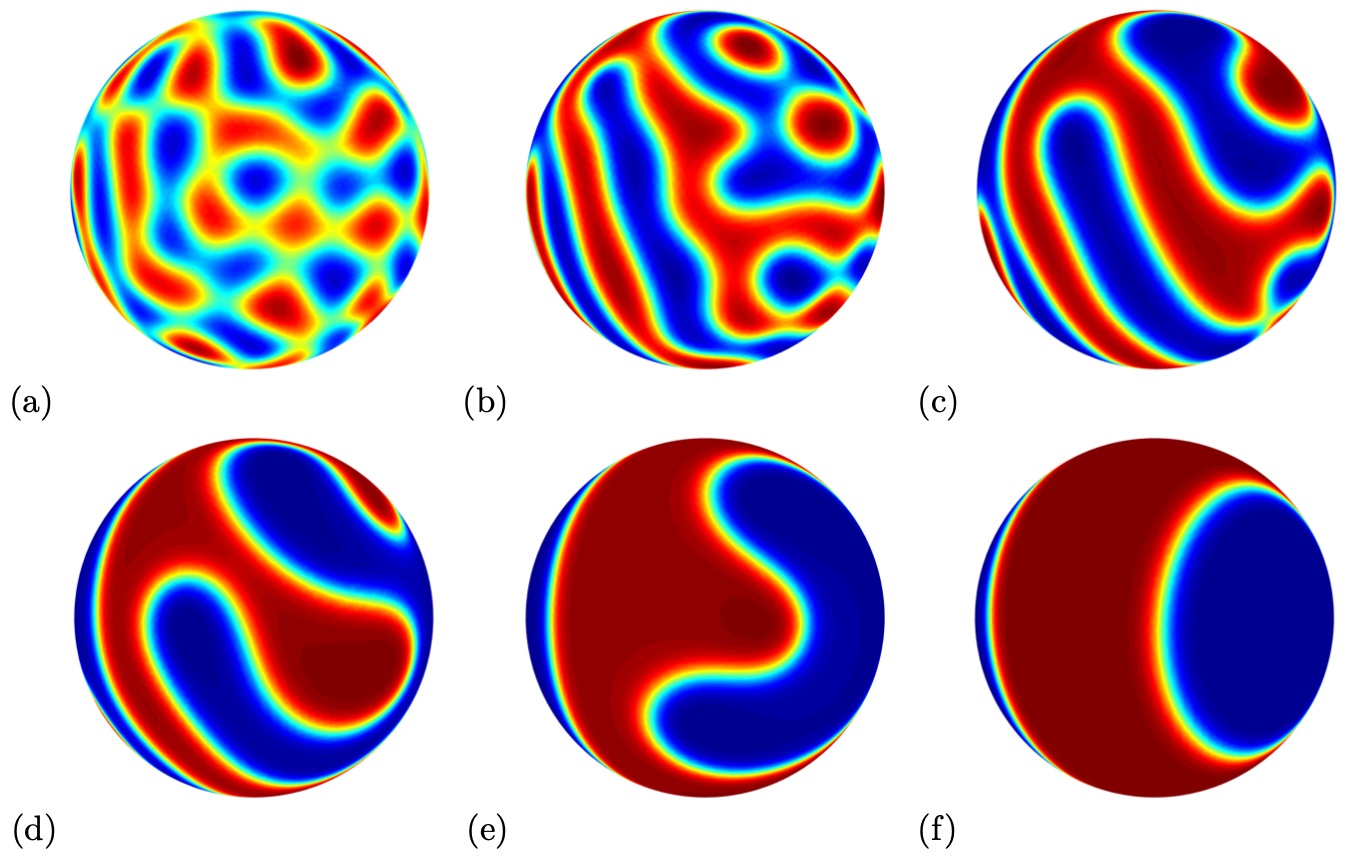}
\caption{(Example \ref{SubSec:RLS-GBPM-Static}) The solution of the Cahn-Hilliard equation on the unit sphere at (a) $t=0.1$, (b) $t=0.2$, (c) $t=0.5$, (d) $t=1$, (e) $t=3$ and (f) $t=9$.}
\label{Fig:CahnSphere}
\end{figure}

\begin{figure}[!ht]
\includegraphics[width=0.95\textwidth]{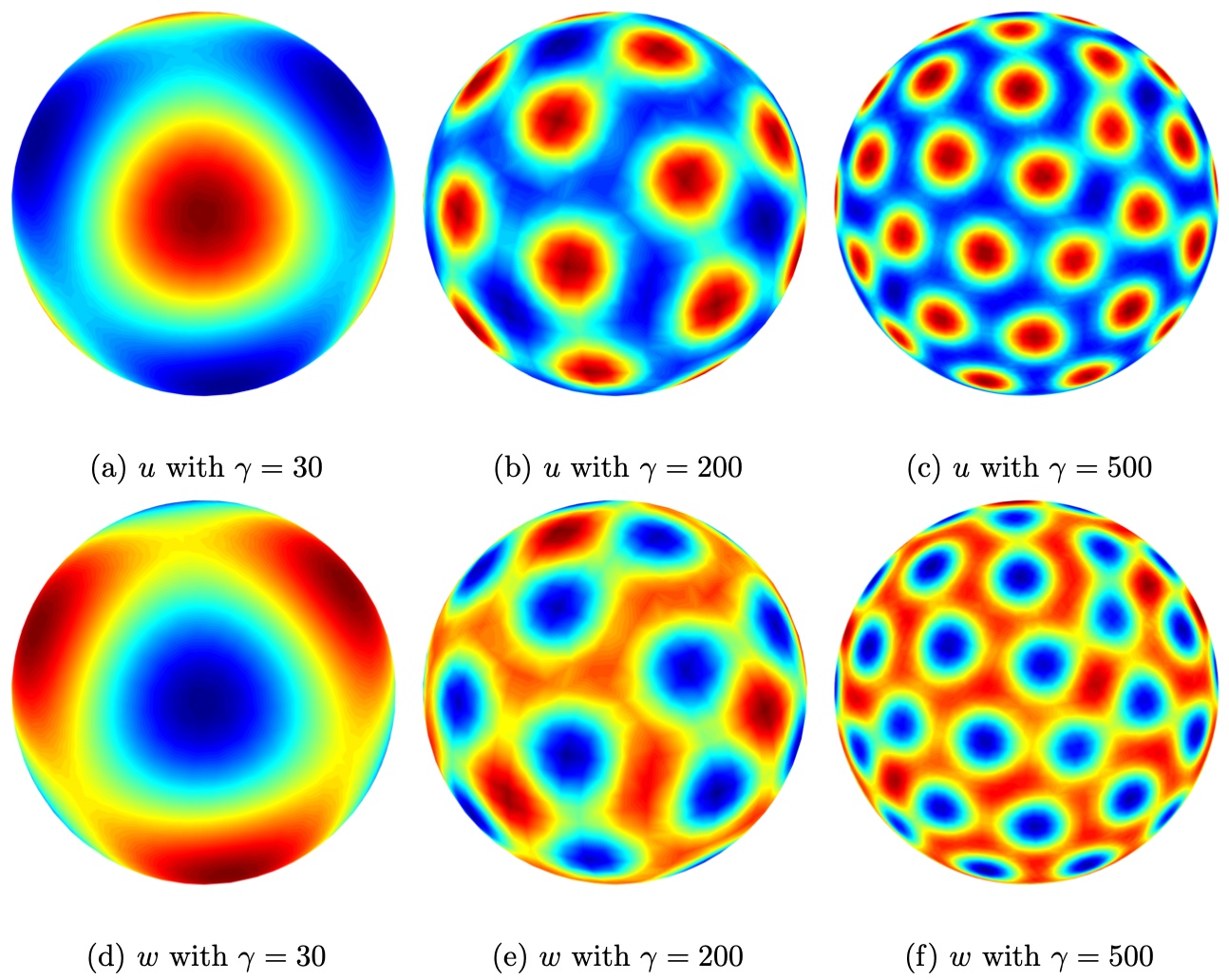}
\caption{(Example \ref{SubSec:RLS-GBPM-Static}) The Turing pattern in the concentrations of (a-c) $u$ and (d-f) $w$ on the unit sphere at $t = 10$ with (from left to right) $\gamma = 30, 200$ and $500$.}
\label{Fig:TuringPattern}
\end{figure}

\subsection{The CLS method for solving PDEs on static surfaces represented by the GBPM}
\label{SubSec:RLS-GBPM-Static}

In this section, we give two numerical examples to show that the proposed CLS method works well in solving PDEs on surfaces sampled by the GBPM representations. The first example is the Cahn-Hilliard equation on the unit sphere. The surface Cahn-Hilliard equation has the form
\begin{eqnarray}
u_t  = \frac{\nu}{Pe}\Delta_{\Gamma}\left[g'(u) \right] - \frac{\nu \, Cn^2}{Pe} \Delta_{\Gamma}^2 u \, .
\label{Eqn.CahnHilliard}
\end{eqnarray}
where $Cn$ is the Cahn-Hilliard number, $Pe$ is the Peclet number and $\nu$ is a diffusion parameter. In this example, we solve the equation on a unit sphere with the parameters in the equation given by $Pe=\nu=1$, $Cn = 0.06$, and the double well potential function $g(u) = u^2(1-u)^2$. Our numerical solution at different times is shown in Figure \ref{Fig:CahnSphere}. The initial condition is a uniform random perturbation with a magnitude of 0.01, with an average concentration of $u$ of 0.5.

In the second example, we simulate the Turing patterns on a static surface. Consider the following coupled system of equation (\ref{Eqn.AdvectionDiffusion}) given by
\begin{eqnarray}
\frac{Du}{Dt} +u\nabla_{\Gamma} \cdot \mathbf{v} - \mathcal{D}_1\Delta_{\Gamma} u &=& f_1(u, w)  \nonumber \\
\frac{Dw}{Dt} +w\nabla_{\Gamma} \cdot \mathbf{v} - \mathcal{D}_2\Delta_{\Gamma} w &=& f_2(u, w) \, ,
\label{Eqn.AdvectionDiffusionSystem}
\end{eqnarray}
with 
\begin{eqnarray}
f_1(u,w) = \gamma(a-u+u^2w) \, \mbox{ and } \, f_2(u,w) = \gamma(b-u^2w) \, .
\label{Eqn.ActivatorDepletedKinetics}
\end{eqnarray}
The parameters $\gamma$, $a$, and $b$ are all positive constants. In this example, we follow \cite{barellmad11} and use the following parameters for the kinetics function in (\ref{Eqn.ActivatorDepletedKinetics}): $a = 0.1$, $b=0.9$, $\mathcal{D}_1 = 1$ and $\mathcal{D}_2 = 10$.  All three examples use the same initial condition given by $(u,w) = (1,0.9)$ with a uniform random perturbation of magnitude 0.01. The parameter $\gamma$ is given in different scales. The grid size in the GBPM representation is $\Delta x= 0.1$, and the stopping time is $t=10$.

Figure \ref{Fig:TuringPattern} shows our numerical solutions corresponding to various $\gamma$ at the fixed final time $t=10$. Even though the exact solution to these examples is not available, these results are similar to those obtained in various articles, including \cite{chagangra01,madmai07,barellmad11}. Numerically, the CLS method is flexible and matches the GBPM well in solving PDEs on a static surface. 


\begin{figure}[!ht]
\centering
\includegraphics[width=0.95\textwidth]{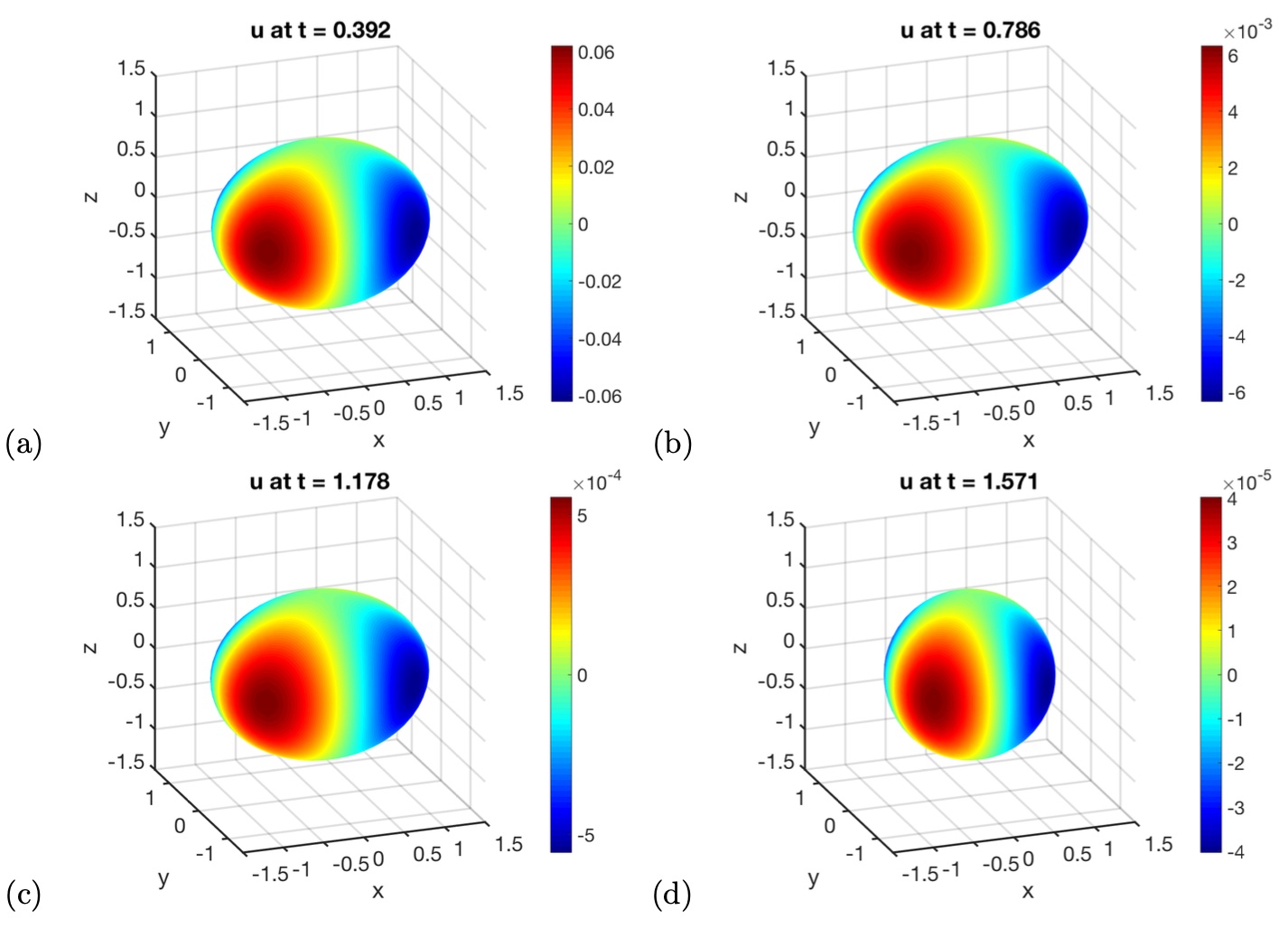}
\caption{(Example \ref{Ex:Convergence}) The numerical solution of the advection-diffusion equation with the source term on an oscillating ellipsoid at (a) $t = \pi/8$, (b) $t = \pi/4$, (c) $t =  3\pi/8$ and (d) $t = \pi/2$.}
\label{Fig:OScillatingSphere}
\end{figure}

\begin{figure}[!ht]
\centering
\includegraphics[width=0.95\textwidth]{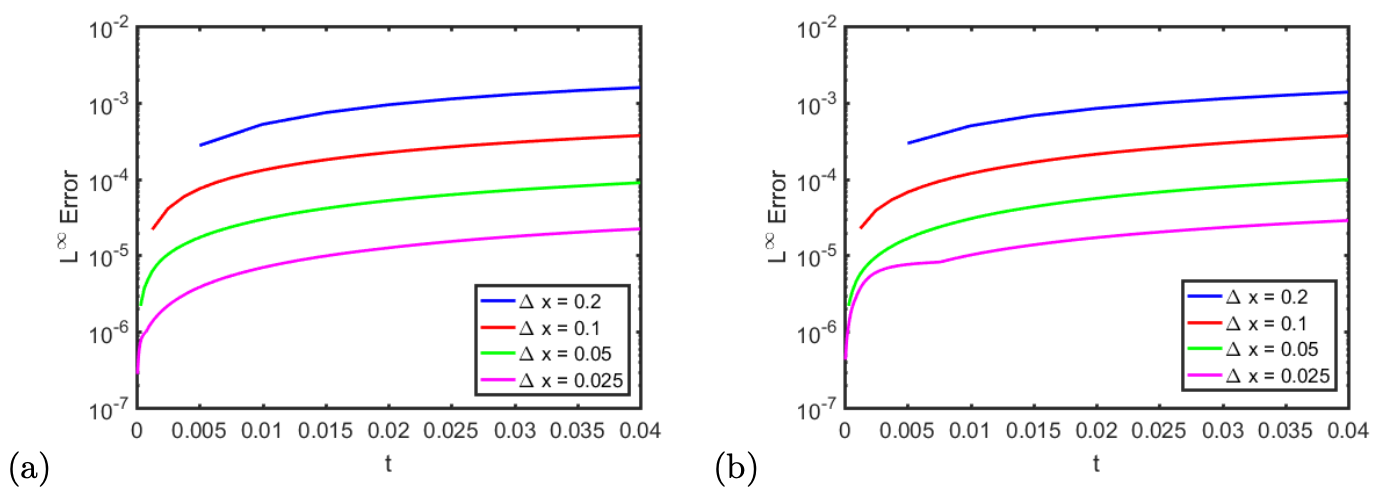}
\caption{(Example \ref{Ex:Convergence}) The $L^{\infty}$-errors in the numerical solutions of the advection-diffusion model on an oscillating ellipsoid over time computed on various meshes. (a) CLS-Polynomial, (b) CLS-RBF.}
\label{Fig:OScillatingSphereConv}
\end{figure}

\begin{table}
\centering
\begin{tabular}{|c|cc|cc|cc|cc|}
  \toprule	
  $\Delta x$ &$ t=0.01$ & order & $t=0.02$ & order & $t=0.03$ & order & $t=0.04$ & order\\
  \hline
  $0.2$ & $5.186\times 10^{-4}$ & - & $9.364\times 10^{-4}$ & - & $1.283\times 10^{-3}$ & - & $1.575\times 10^{-3}$ & -\\
  $0.1$ & $1.299\times 10^{-4}$ & $2.000$ & $2.222\times 10^{-4}$ & $2.108$ & $3.012\times 10^{-4}$ & $2.130$ & $3.700\times 10^{-4}$ & $2.130$\\
  $0.05$ & $2.956\times 10^{-5}$ & $2.197$ & $5.201\times 10^{-5}$ & $2.136$ & $7.169\times 10^{-5}$ & $2.101$ & $8.927\times 10^{-5}$ & $2.072$\\
  $0.025$ & $6.909\times 10^{-6}$ & $2.139$ & $1.253\times 10^{-5}$ & $2.076$ & $1.769\times 10^{-5}$ & $2.026$ & $2.221\times 10^{-5}$ & $2.010$\\
  \bottomrule
\end{tabular}
\caption{(Example \ref{Ex:Convergence}, CLS-Polynomial) $L_{\infty}$-errors and some estimated orders of convergence at different times.}
\label{TB:Convergence}
\end{table}

\begin{table}
\centering
\begin{tabular}{|c|cc|cc|cc|cc|}
  \toprule	
  $\Delta x$ &$ t=0.01$ & order & $t=0.02$ & order & $t=0.03$ & order & $t=0.04$ & order\\
  \hline
  $0.2$ & $4.957\times 10^{-4}$ & - & $8.382\times 10^{-4}$ & - & $1.122\times 10^{-3}$ & - & $1.372\times 10^{-3}$ & -\\
  $0.1$ & $1.177\times 10^{-4}$ & $2.106$ & $2.109\times 10^{-4}$ & $1.987$ & $2.939\times 10^{-4}$ & $1.909$ & $3.671\times 10^{-4}$ & $1.868$\\
  $0.05$ & $3.017\times 10^{-5}$ & $1.950$ & $5.563\times 10^{-5}$ & $1.896$ & $7.824\times 10^{-5}$ & $1.878$ & $9.980\times 10^{-5}$ & $1.873$\\
  $0.025$ & $1.007\times 10^{-5}$ & $1.497$ & $1.706\times 10^{-5}$ & $1.630$ & $2.313\times 10^{-5}$ & $1.691$ & $2.857\times 10^{-5}$ & $1.715$\\
  \bottomrule
\end{tabular}
\caption{(Example \ref{Ex:Convergence}, CLS-RBF) $L_{\infty}$-errors and some estimated orders of convergence at different times.}
\label{TB:ConvergenceRBF}
\end{table}

\subsection{The advection-diffusion equation with a source term on an oscillating ellipsoid}
\label{Ex:Convergence}
This three-dimensional example is taken from \cite{dziell07,ellstisty10,petruu16} to test the accuracy of a PDE solver on moving surfaces. The time-dependent oscillating ellipsoid has the following exact solution

$$
\Gamma(t) = \left\lbrace \x = (x,y,z) \in \mathbb{R}^3\left| \frac{x^2}{a(t)}+y^2+z^2 = 1\right. \right\rbrace \, ,
$$
where $a(t) = 1+\sin 2t$. The surface $\Gamma$ can be explicitly parametrized as

$$
\Gamma(t,\theta,\phi) := \left(\sqrt{a(t)}\cos\theta\cos\phi, \sin\theta\cos\phi,  \sin\phi \right) \, .
$$
The associated velocity field is given by $\mathbf{v}(\x,t) = \frac{a'(t)}{a(t)} (x,0,0)^T$. Therefore, the equation involves both the advection along the surface and the deformation by the motion in the normal direction. We choose the exact continuous solution $u(\x,t) = e^{-6t}xy$ so that the corresponding source term on the right-hand side of the equation is given by
\begin{eqnarray*}
f(\x,t) = \left\{
-6+\frac{a'(t)}{a(t)}
\left( 1-\frac{x_1^2}{2N} \right)+\frac{1+5a(t)+2a^2(t)}{N} - \frac{1-a(t)}{N^2} \left[x_1^2+a^3(t) \left(x_2^2+x_3^2\right) \right] \right\}  u(\x,t)
\end{eqnarray*}
with $N = x^2+a^2(t)(y^2+z^2)$. In this simulation, we choose the mesh size $\Delta x = 0.1$ and a fixed time step size $\Delta t = 0.1 \, (\Delta x)^2$. We plot our solution at different time $t = \pi/8, \pi/4, 3\pi/8$ and $\pi/2$ in Figure \ref{Fig:OScillatingSphere}. Figure \ref{Fig:OScillatingSphereConv} gives the $L^{\infty}$-error in the numerical solution computed using the CLS with different mesh sizes at different times. Table \ref{TB:Convergence} shows the errors in the computed solution at different times $t = 0.01, 0.02, 0.03$ and $0.04$ computed on a different mesh. We see that our numerical approach is roughly second-order accurate.

Concerning the computational efficiency, the CPU time for one marching step on the mesh with $\Delta x = 0.1$ is $3.1$ seconds, including $1.5$ seconds spent on the CLS reconstruction and $1.6$ seconds on the GBPM computations. This computational time is recorded under a 4-core laptop in \textsf{MATLAB} with a parallel modulo implemented. 


\begin{figure}[!ht]
\centering
\includegraphics[width=0.95\textwidth]{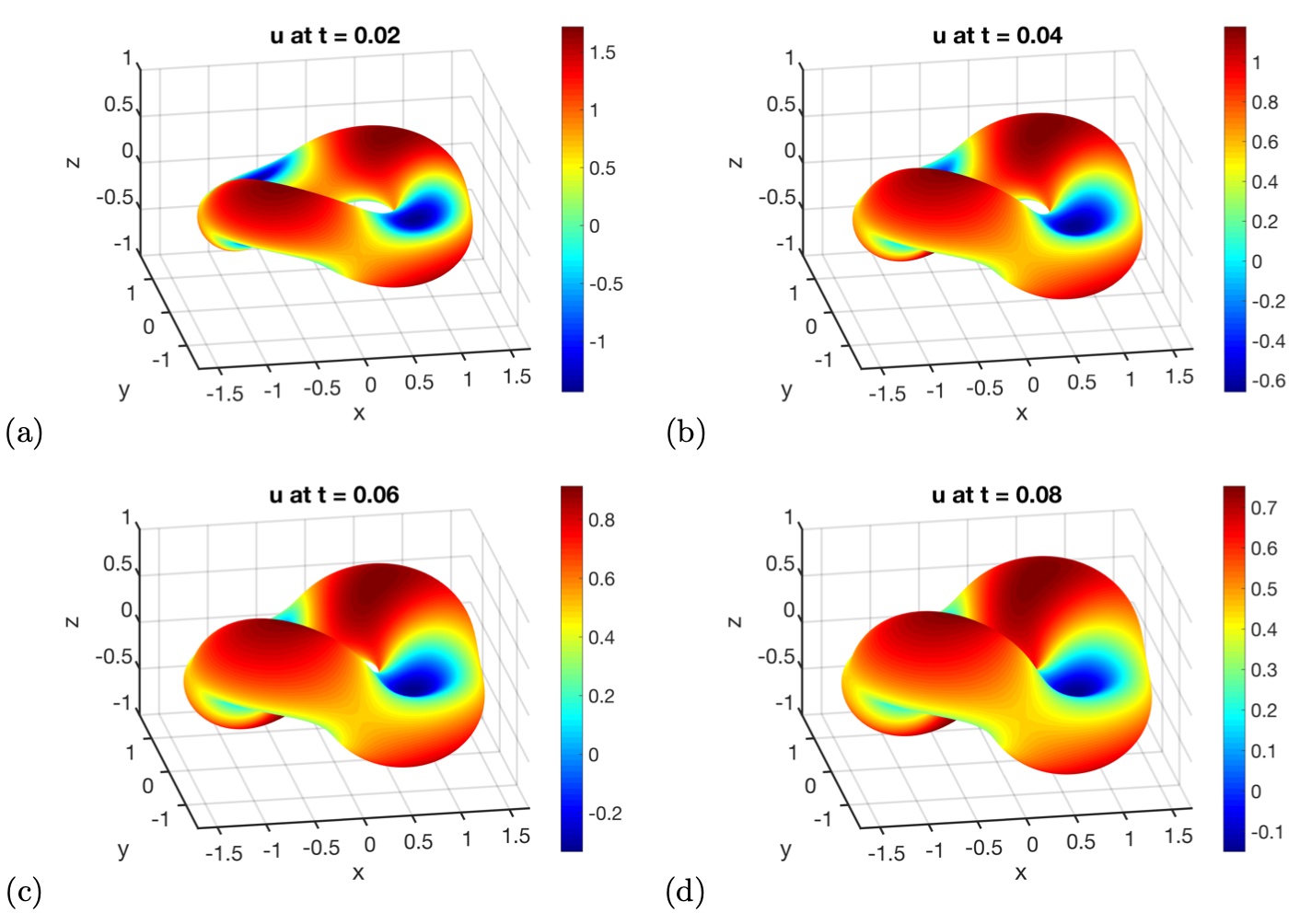}
\caption{(Example \ref{Ex:CoupledEvolvingTorus}) The numerical solution $u$ at (a) $t = 0.02$, (b) $t = 0.04$, (c) $t = 0.06$ and (d) $t = 0.08$.}
\label{Fig:CoupledTorus}
\end{figure}

\subsection{A strongly coupled flow on an evolving torus}
\label{Ex:CoupledEvolvingTorus}
In this case, the interface velocity is strongly coupled with the solution of the PDE. This example is taken from \cite{petruu16}, where we evolve an initial torus 

$$
\left(1-\sqrt{x^2+y^2}\right)^2+z^2 = 0.3^2
$$ 
according to the given velocity $\mathbf{v} = (0.1\kappa + 5u)\n$ where $\kappa$ is the mean curvature, $\n$ is the unit normal vector, and $u$ is the solution to the PDE 
\begin{eqnarray*}
u_t + V\frac{\partial u}{\partial \n} -\kappa V u + \Delta_{\Gamma} u = 0 \text{\;  on\; } \Gamma \, .
\end{eqnarray*}
The initial value of $u$ is given by $u(\x,t) = 1+ 20xyz$. Figure \ref{Fig:CoupledTorus} shows our computed solution of $u$ on the evolving surface at different times from $0.02$ to $0.08$ in an increment of $0.02$. Even though the exact solution to this problem is unavailable, these numerical solutions are similar to those presented in \cite{petruu16}.


\begin{figure}[!ht]
\centering
\includegraphics[width=0.95\textwidth]{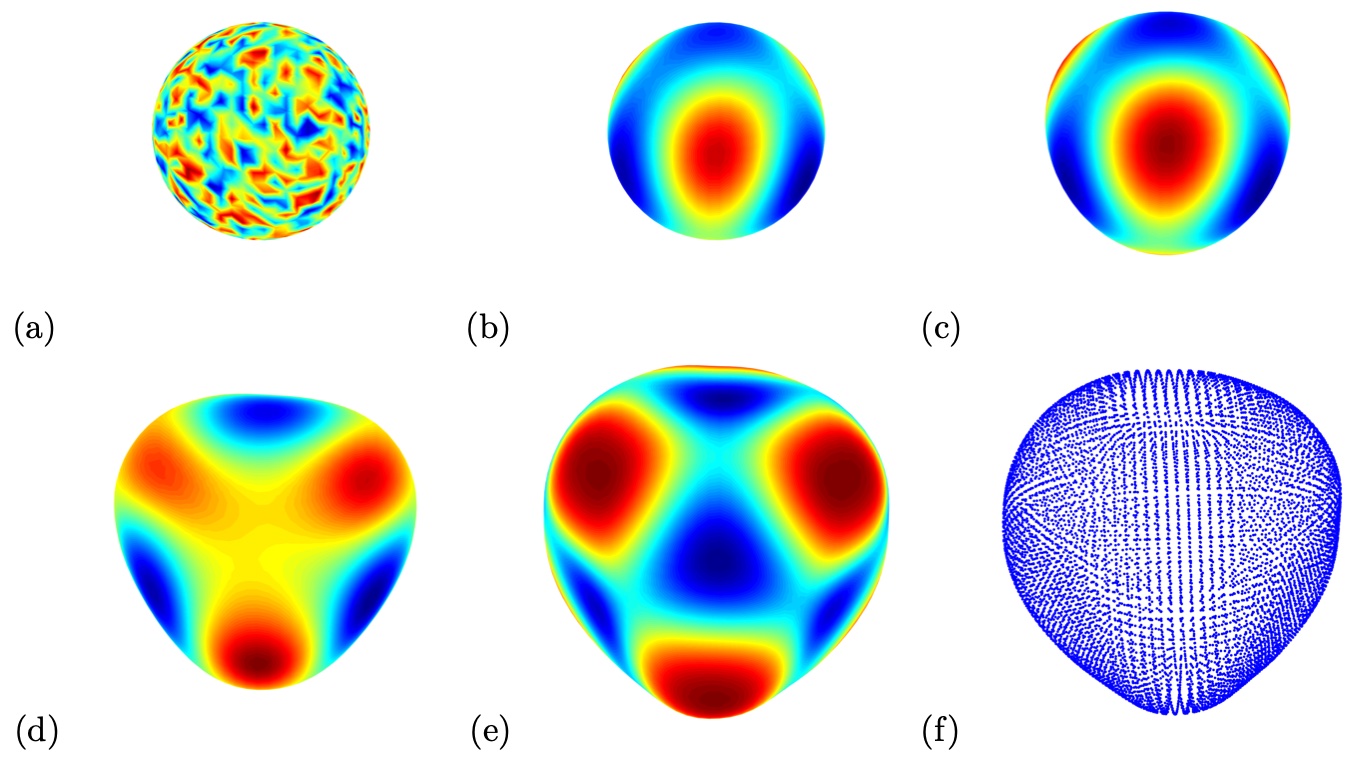}
\caption{(Example \ref{Ex:GrowthTumor}) The tumor surface and the concentration of the growth factor $u$ at (a) $t=0$, (b) $t=5$, (c) $t=6$, (d) $t=8$ and (e) $t=10$. (f) The corresponding foot points at final time $t=10$. The parameters used in this example are $\gamma=30$, $\mathcal{D}_1=1$, $\mathcal{D}_2=10$, $a=0.1$, $b=0.9$, $\bar{t}=5$, $\delta=0.1$, $\epsilon=0.01$ and the time marching step is 0.0002. The number of footpoints grows from $3810$ to $9730$ at the final time.}
\label{Fig:GrowthTumor1}
\end{figure}

\begin{figure}[!ht]
\centering
\includegraphics[width=0.95\textwidth]{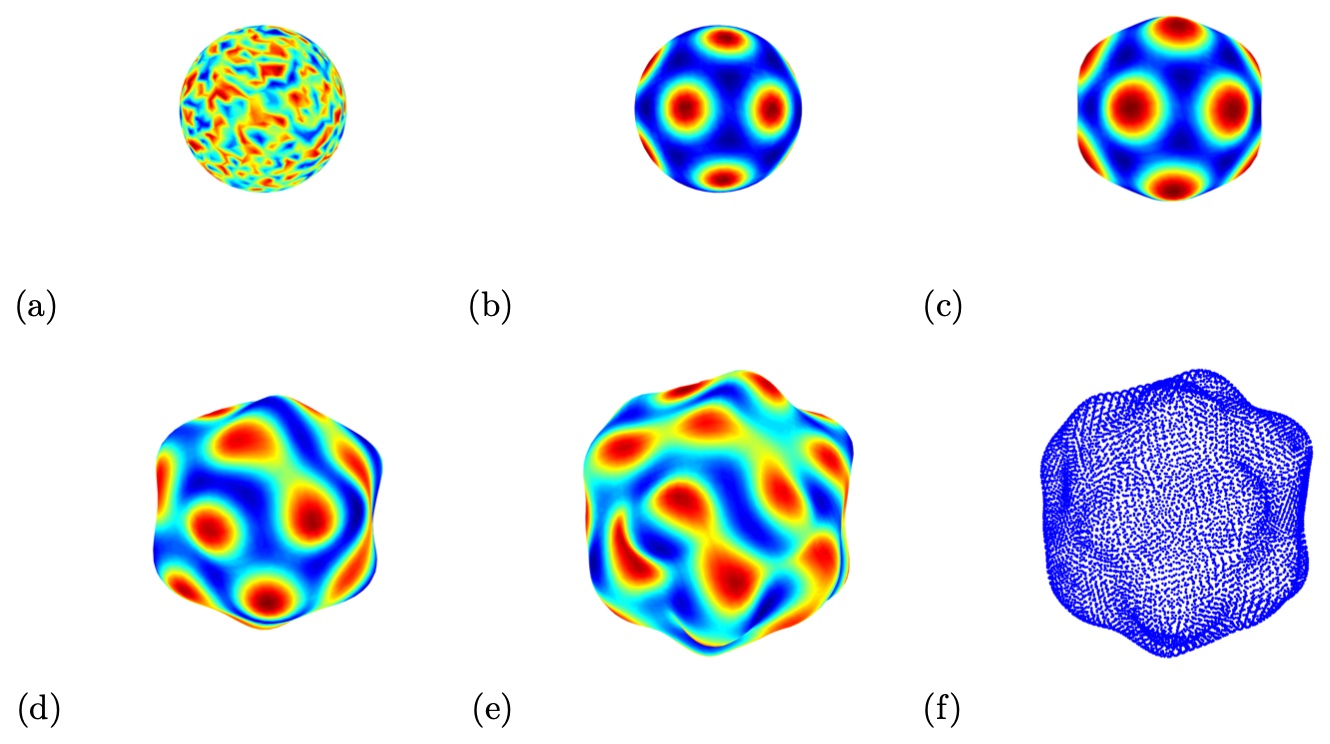}
\caption{(Example \ref{Ex:GrowthTumor}) The tumor surface and the concentration of the growth factor $u$ at (a) $t=0$, (b) $t=5$, (c) $t=6$, (d) $t=8$ and (e) $t=10$. (f) The corresponding foot points at final time $t=10$. The parameters used in this example are $\gamma=100$, $\mathcal{D}_1=1$, $\mathcal{D}_2=10$, $a=0.1$, $b=0.9$, $\bar{t}=5$, $\delta=0.1$, $\epsilon=0.01$ and the time marching step is 0.0002. The number of footpoints grows from $3810$ to $10745$ at the final time.}
\label{Fig:GrowthTumor2}
\end{figure}

\subsection{Solid tumors growth} 
\label{Ex:GrowthTumor}
In this example, we compute the reaction-diffusion system (\ref{Eqn.AdvectionDiffusionSystem}) with an activator-depleted kinetics in the form of (\ref{Eqn.ActivatorDepletedKinetics}). The surface growth law is assumed to be in the form 
\begin{eqnarray*}
\mathbf{v} = (-\epsilon \kappa + \delta u)\n, \;\; t\geq \bar{t}
\end{eqnarray*}
where $\delta$ is the growth rate and the parameter $\epsilon$ models the surface tension of the tumor surface. In this system, the solution of $u$ denotes the growth-promoting factor, and $w$ is the corresponding growth-inhibiting factor \cite{cragafmai99}. We follow the discussions in \cite{chagangra01,barellmad11,ellsty12} by first starting the numerical simulations from a stationary sphere in a pre-patterning stage (controlled by the parameter $\bar{t}$). After the pre-patterning stage, the surface grows under the given motion law.

Figures \ref{Fig:GrowthTumor1} and \ref{Fig:GrowthTumor2} show the evolution of the tumor surface together with the growth factor concentration $u$ with two different values of $\gamma$'s. In the pre-pattern stage $t\in[0,\bar{t}]$, the behavior of the solution $u$ is the same as when we simulate the Turing patterns in the previous section. After the pre-pattern stage, the surface evolves in the normal direction. Those regions with a larger growth concentration $u$ evolve much faster than others. This leads to some interesting configurations of the shape.


\begin{figure}[!ht]
\includegraphics[width=0.95\textwidth]{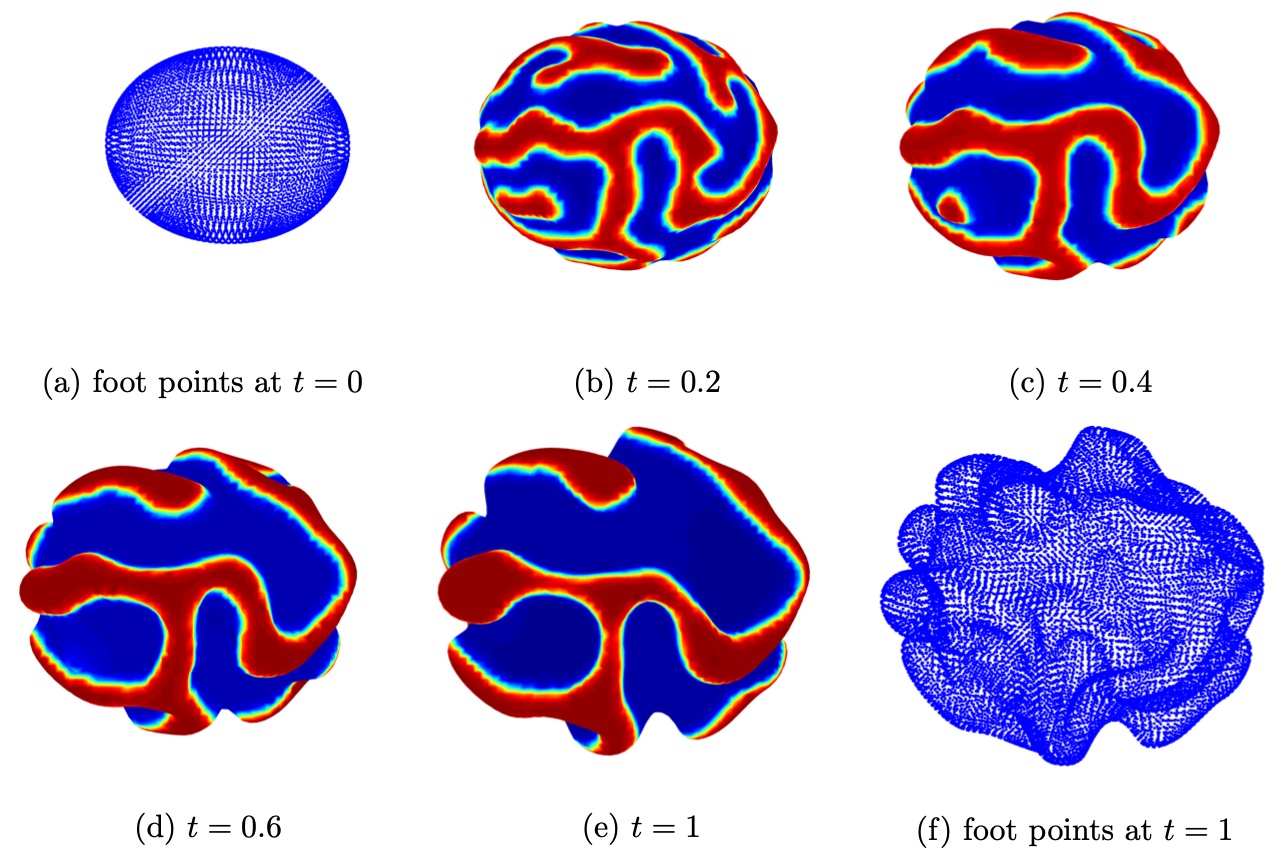}
\caption{(Example \ref{Ex:CahnHilliardEllipsoid}) The evolution of the interface and the solution $u$ to the Cahn-Hilliard equation with $\gamma = 100$ at the time $t=0.2$, 0.4, 0.6 and 1.0. The number of footpoints grows from (a) $11312$ at $t=0$ to (f) $22062$ at the final time $t=1$.}
\label{Fig:CahnHilliardEllipsoid}
\end{figure}

\subsection{The Cahn-Hilliard equation on an ellipsoid}
\label{Ex:CahnHilliardEllipsoid}
This example considers the homogenous Cahn-Hilliard equation (\ref{Eqn.CahnHilliard}) on an evolving ellipsoid. An initial ellipsoid centered at the origin with axis $r_x = 1$ and $r_y = r_z = 0.8$ is evolved in the normal direction according to the velocity $\mathbf{v} = (0.01\kappa + 0.4u)\n$, where $\kappa$ is the mean curvature and $\n$ is the unit normal direction. We choose the Peclet number $Pr=1$, the diffusion parameter $\nu=1$, $Cn = 0.03$ and the double well potential function $g(u) = u^2(1-u)^2$. The initial value of the function $u$ is defined as a $0.01$ Gaussian noise around a constant number $0.5$. Figure \ref{Fig:CahnHilliardEllipsoid} shows the numerical results computed on the mesh with $\Delta x = 0.05$ and the time marching step $\Delta t =0.0001$.

\section{Conclusion}
This paper proposes a new CLS method for solving PDEs on evolving surfaces. The method modifies the GBPM developed in \cite{leuzha0801,leuzha0802,leulowzha11} by enforcing an extra constraint in the local least squares approximation. Unlike a recent approach in \cite{wanleuzha17} where we introduce extra weighting to the target point through a virtual grid, the additional condition in the current approach requires the least squares reconstruction passing through the data point at the target location. This significantly improves the computational stability and accuracy. We have shown that this CLS approach fits extremely well with the GBPM for interface representation. The coupling of these two methods provides a computationally efficient and robust method for solving PDEs on evolving surfaces.

\section*{Acknowledgment}
Leung's work was partly supported by the Hong Kong RGC grants 16303114 and 16309316.

\bibliographystyle{plain}
\bibliography{syleung,mybib}

\end{document}